\documentclass[a4paper, notitlepage, 11pt]{article}

\usepackage{ucs}                
\usepackage[utf8x]{inputenc} 	

\usepackage{color}
\usepackage{lscape}
\usepackage{amsmath}
\usepackage{amssymb}
\usepackage{amsbsy}
\usepackage{amsthm}
\usepackage{wasysym}
\usepackage{varioref}
\usepackage{tocbibind}          
\usepackage{float}
\usepackage[square,numbers]{natbib}
\usepackage{url}
\usepackage{authblk}

\usepackage[british]{babel}

\usepackage{refdef}
\usepackage{nfontdef}

\newcommand{\ignore}[1]{}

\newcommand{\vrho}{\varrho}
\newcommand{\vepsilon}{\varepsilon}
\newcommand{\vsigma}{\varsigma}

\newcommand{\vphi}{\varphi}

\newcommand{\vOmega}{\varOmega}

\newcommand{\vek}[1]{\mathchoice{\displaystyle\boldsymbol{#1}}
{\textstyle\boldsymbol{#1}}{\scriptstyle\boldsymbol{#1}}
{\scriptscriptstyle\boldsymbol{#1}}}
\newcommand{\mat}[1]{\mathchoice{\displaystyle\mathbf{#1}}
{\textstyle\mathbf{#1}}{\scriptstyle\mathbf{#1}}
{\scriptscriptstyle\mathbf{#1}}}
\newcommand{\opb}[1]{\vek{{\mathsf{#1}}}}

\newcommand{\EXP}[1]{\mathbb{E}\left(#1\right)}
 
\newcommand{\dlangle}{\langle\negthinspace\langle}
\newcommand{\drangle}{\rangle\negthinspace\rangle}

\newcommand{\divg}{\mathop{\mathrm{div}}\nolimits}

\newcommand{\dd}{\partial}
\newcommand{\di}{\mathrm{d}}

\newcommand{\KL}{Karhunen-Lo\`eve}
\newcommand{\ip}[2]{\langle #1, #2 \rangle}
\newcommand{\bkt}[2]{\langle #1 | #2 \rangle}

\newcommand{\ns}[1]{| #1 |}
\newcommand{\nd}[1]{\| #1 \|}

\definecolor{myred}{rgb}{1, 0.2, 0.2}

\newcommand{\thetitle}{Parameter~Identification~in~a
                       Probabilistic~Setting}

\newcommand{\theauthorA}{Bojana~V.~Rosi\'c}
\newcommand{\theauthorB}{Anna~Ku\v{c}erov\'a}
\newcommand{\theauthorC}{Jan~S\'ykora}
\newcommand{\theauthorD}{Oliver~Pajonk}
\newcommand{\theauthorE}{Alexander~Litvinenko}
\newcommand{\theauthorF}{Hermann~G.~Matthies}

\newcommand{\theauthor}{
\theauthorA, \theauthorB, 
\theauthorC, \theauthorD,
\theauthorE, \theauthorF }

\newcommand{\thesubject}{Parameter identification, Bayesian updating}
\newcommand{\theclassification}{{{\bf (MSC2010) 62F15}, 65N21, 62P30, 60H15, 60H25, 74G75, 80A23, 74C05}\\
{{\bf (PACS2010) 46.65.+g}, 46.35.+z, 44.10.+i}\\
{{\bf (ACM1998)} G.1.8, G.3, J.2}}
\newcommand{\thekeywords}{parameter identification, Bayesian update,
                         linear Bayes, Kalman filter, polynomial chaos}
\newcommand{\textdate}{December 2011}

\newcommand{\thebib}{./}

\newcommand{\pdflinkcolor}{blue}

\usepackage{ifpdf}

\ifpdf
  \usepackage[pdftex,
       colorlinks=true,		
       bookmarksnumbered,	
       bookmarksopen,		
       pdfstartview=FitH,	
       linkcolor=\pdflinkcolor,	
       citecolor=\pdflinkcolor,	
       urlcolor=\pdflinkcolor,	
       filecolor=\pdflinkcolor,	
       hyperindex,		
       plainpages=false,	
       breaklinks=true,		
       pdfpagelabels, 
       hypertexnames=true,
     ]{hyperref}
     \usepackage[dvipdf]{graphicx}
	 \usepackage[pdftex]{thumbpdf}
\else
  \usepackage[
    colorlinks=true,
    linkcolor=\pdflinkcolor,
    hyperindex,
    plainpages=false
  ]{hyperref}
  \usepackage[dvips]{graphicx}
\fi

\hypersetup{a5paper,
            colorlinks,
            linkcolor=\pdflinkcolor,
            citecolor=\pdflinkcolor,
            urlcolor=\pdflinkcolor,
            pdftitle={\thetitle},
            pdfauthor={\theauthor},
            pdfsubject={\thesubject},
            pdfkeywords={\thekeywords},
            pdfproducer={\theauthor~using \LaTeXe}}

\newcommand{\HRule}{\rule{\linewidth}{0.8mm}}

\begin{document}

\title{\thetitle}
\author[a,b]{\theauthorA}
\author[c]{\theauthorB}
\author[c]{\theauthorC}
\author[d,a]{\theauthorD}
\author[a]{\theauthorE}
\author[a]{\theauthorF
}

\date{\textdate}

\makeatletter
\affil[a]{Technische Universit\"at Braunschweig\authorcr
\texttt{\href{mailto:wire@tu-bs.de?subject=\thetitle}{wire@tu-bs.de}}}
\affil[b]{University of Kragujevac}
\affil[c]{Czech Technical University in Prague}
\affil[d]{SPT Group Hamburg}
\makeatother

\ignore{


\setcounter{page}{0}
\thispagestyle{empty}
\cleardoublepage

\begin{titlepage}
	\begin{center}
	{\Large \textsc{Bojana~V.~Rosi\'c, Anna~Ku\v{c}erov\hspace{-1.5pt}\'a,
     Jan~S\'ykora, Oliver~Pajonk, Alexander~Litvinenko, Hermann~G.~Matthies}}
	\end{center}
	
   \vspace*{\stretch{1}}
	\begin{center}
	\HRule
	\end{center}
	\begin{flushright}
		\huge\scshape 
		 \thetitle
		 \\[5mm]
	\end{flushright}
	\begin{center}
	\HRule
	\end{center}
	
	\vspace*{\stretch{3}}
	

	\vspace*{\stretch{3}}

	\begin{tabular}{l|r}
\begin{minipage}{2.5cm}
 \includegraphics[width=2.5cm]{common/tu-bs-color-logo}
\end{minipage}&
\hspace{0.1cm}
\begin{minipage}{10 cm}
{\large \textsc{Informatikbericht Nr. \thereport}\\[3mm]
\textsc{Institute of Scientific Computing}\\
\textsc{Carl-Friedrich-Gauß-Fakult\"at}\\
\textsc{Technische Universit\"at Braunschweig}}
\\\
\\\
Brunswick, Germany\\
\end{minipage}
\end{tabular}
\end{titlepage}

\newpage

\thispagestyle{empty}
\vspace*{\stretch{2}}

\begin{flushleft}
\begin{tabular}{ll}
\makeatletter
This document was created \textdate{} using \LaTeXe. \\[1cm]
\makeatother
\end{tabular}

\begin{tabular}{ll}
\begin{minipage}{6cm}
Institute of Scientific Computing\\ 
Technische Universit\"at Braunschweig\\
Hans-Sommer-Stra\ss{}e 65\\
D-38106 Braunschweig, Germany\\

\texttt{url: \url{www.wire.tu-bs.de}}\\
\makeatletter
\texttt{mail: \href{mailto:wire@tu-bs.de?subject=\thetitle}{wire@tu-bs.de}}
\makeatother
\end{minipage}
&
\begin{minipage}{2.5cm}
\vspace{-0.5cm}
\includegraphics[scale=0.34]{common/logo_wire}

\end{minipage}
\end{tabular}

\vspace*{\stretch{1}}

Copyright \copyright{} by \theauthor{}\\[5mm]
\end{flushleft}

This work is subject to copyright. All rights are reserved, whether the whole or part of the material is concerned, specifically the rights of translation, reprinting, reuse of illustrations, recitation, broadcasting, reproduction on microfilm or in any other way, and storage in data banks. Duplication of this publication or parts thereof is permitted in connection with reviews or scholarly analysis. Permission for use must always be obtained from the copyright holder.\\[5mm]

Alle Rechte vorbehalten, auch das des auszugsweisen Nachdrucks, der auszugsweisen oder vollständigen Wiedergabe (Photographie, Mikroskopie), der Speicherung in Datenverarbeitungsanlagen und das der Übersetzung.


}

\maketitle

\begin{abstract}
Parameter identification problems are formulated in a probabilistic
language, where the randomness reflects the uncertainty about the
knowledge of the true values.  This setting allows conceptually easily
to incorporate new information, e.g.\ through a measurement, by 
connecting it to Bayes's theorem. The unknown quantity is modelled as
a (may be high-dimensional) random variable.  Such a description has two
constituents, the measurable function and the measure.  One group of
methods is identified as updating the measure, the other group changes
the measurable function.  We connect both groups with the relatively
recent methods of functional approximation of stochastic problems, and
introduce especially in combination with the second group of methods
a new procedure which does not need any sampling, hence works completely
deterministically.  It also seems to be the fastest and more reliable when
compared with other methods.  We show by example that it also works
for highly nonlinear non-smooth problems with non-Gaussian measures.

\vspace{5mm}

\vspace{5mm}
{\noindent\textbf{Classification:}\\ \theclassification}

\end{abstract}

\pagenumbering{roman}

\pagenumbering{arabic}

%

\section{Introduction}  \label{S:intro}
In trying to predict the behaviour of physical systems, one is often
confronted with the fact that although one has a mathematical model
of the system which carries some confidence as to its fidelity, some
quantities which characterise the system may only be incompletely
known, or in other words they are uncertain.
See \cite{boulder:2011} for a synopsis on our approach to
such parametric problems.

Here we want to identify these parameters through
observations or measurement of the response of the system.
Such an identification can be approached in different ways.
One way is to measure the difference between observed
and predicted system output and try to find parameters such
that this difference is minimised, this optimisation
approach leads to regularisation procedures \cite{Engl2000}.

Here we take the view that our lack of knowledge or uncertainty
of the actual value of
the parameters can be described in a \emph{Bayesian} way through a
probabilistic model \cite{jaynes03,Tarantola2004}.  The unknown parameter
is then modelled as a random variable---also called the \emph{prior} model---and
additional information on the system through measurement or observation
changes the probabilistic description to the so-called \emph{posterior} model.
The second approach is thus a method
to update the probabilistic description in such a way as to take account
of the additional information.

To be more specific, let us consider the following situation:
we are investigating some physical system which
is modelled by an evolution equation for its state: 
\begin{equation} \label{eq:I}
\frac{\dd}{\dd t}u(t) + A(p;u(t)) = f(p;t),
\end{equation}
where $u(t) \in \C{U}$ describes the state of the system
at time $t \in [0,T]$ lying in a Hilbert space $\C{U}$ (for the sake of
simplicity), $A$ is an
operator modelling the physics of the system, and $f\in\C{U}^*$ is some external
influence (action / excitation / loading).  The model depends on some parameters
$p \in \C{P}$, and by $q$ we denote that component of the parameters $p$
which we are uncertain about and would thus like to identify the actual value.

Now assume that we observe a function of
the state $Y(u(q),q)$, and from this observation we would like to identify
the corresponding $q$.  This is called the \emph{inverse} problem, and
as the mapping $q\mapsto Y(q)$ is usually not invertible, the inverse
problem is \emph{ill-posed}.  Embedding this problem of finding the best
$q$ in a larger class by modelling our knowledge
about it with the help of probability theory, then in a Bayesian manner the
task becomes to estimate conditional expectations, e.g.\ see 
\cite{jaynes03,Tarantola2004}
and the references therein.  The problem now is
\emph{well-posed}, but at the price of `only' obtaining probability
distributions on the possible values of $q$, which now is modelled
as a $\C{Q}$-valued random variable (RV).  Predicting what the measurement
$Y(q)$ should be from some assumed $q$ is computing the \emph{forward}
problem.  The inverse problem is then approached by comparing the
forecast from the forward problem with the actual information.

Since the parameters of the model to be estimated are uncertain, all relevant
information may be obtained via their stochastic description.
In order to extract information from the posterior most estimates take
the form of expectations w.r.t.\ the posterior.
These expectations---mathematically integrals, numerically to be evaluated
by some quadrature rule---may be computed via asymptotic,
deterministic or sampling methods.  In our review of current work we
follow our recent reports \cite{bvrAlOpHgm11,opBvrAlHgm11}.

One often used technique is a 
Markov chain Monte Carlo (MCMC) method \cite{Madras-Fields:2002,Gamerman06},
constructed such
that the asymptotic distribution of the Markov chain is the Bayesian
posterior distribution.
This can be then sampled by letting the Markov chain run for a sufficiently
long time, although the samples are not independent in this case.
With the intention of accelerating the MCMC method some authors
\cite{Marzouk2009a,Kucherova10,Pence2010,kucerovaSykBvrHgm:2011}
have introduced stochastic spectral methods into the
computation.  Expanding the prior random process into a polynomial chaos (PCE) 
or a \KL{} expansion (KLE) (e.g.\ \cite{matthies6}), the inverse problem
becomes an inference on the weights of the \KL{} modes.
Pence et al.\ \cite{Pence2010} combine polynomial chaos theory with maximum likelihood estimation, where the parameter estimates are calculated in a recursive
or iterative manner.  Christen and Fox \cite{Christen2005} have applied a local
linearisation of the forward model to improve the acceptance probability of proposed moves, while in \cite{Balakrishnan2003,Ma2009,Li2009,zengZhang:2010}
collocation methods are employed as a more efficient sampling technique.

The approaches mentioned above require a large number of samples in order
to obtain satisfactory results.  While showing some results in this direction,
the main idea here is to do the Bayesian update directly on the
polynomial chaos expansion (PCE) without any sampling 
\cite{opBvrAlHgm11,bvrAlOpHgm11,boulder:2011}.
This idea has appeared independently in \cite{Blanchard2010a}
in a simpler context, whereas
in \cite{saadGhn:2009} it appears as a variant of the Kalman filter
(e.g.\ \cite{Kalman}).
A PCE for a push-forward of the posterior
measure is constructed in \cite{moselhyYMarz:2011}.

From this short overview it becomes apparent that
the update may be seen abstractly
in two different ways.
Regarding the uncertain parameters 
\begin{equation}  \label{eq:RVq}
q: \vOmega \to \C{Q} \text{  as a RV on a probability space   }
  (\vOmega, \F{A}, \D{P})
\end{equation}
where the set of elementary events is $\vOmega$, $\F{A}$ a $\sigma$-algebra of
events, and $\D{P}$ a probability measure, one set of methods performs
the update by changing the probability measure $\D{P}$ and leaving the
mapping $q(\omega)$ as it is, whereas the other set of methods leaves the
probability measure unchanged and updates the function $q(\omega)$.
In any case, the push forward measure $q_* \D{P}$ on $\C{Q}$ is changed
from prior to posterior.

The organisation of the paper is as follows: in \refS{bayes} we review
the Bayesian update and recall the link between the conditional measure
and conditional expectation.  This allows to recover a Kalman filter like
update for the RV via the conditional expectation.  In the following
\refS{num-real} different ways of numerically computing the Bayesian
update are indicated.  The resolution of the forward problem in the context
of stochastic discretisation procedures is outlined in \refS{forward}.
The numerical examples for such parameter identification procedures,
contained in \refS{num-xmpls}, are representative of some linear and
nonlinear problems in structural and continuum mechanics.

%
%
%
%
%


%

\section{Bayesian Updating} \label{S:bayes}
In the setting of \refeq{eq:I} let us pose the following problem: Some
components---let us denote these by $q$---of the parameters
$p \in \C{P}$ are uncertain.  To be more specific, assume that 
$q \in \C{Q}$ are elements of some vector space.
By making observations $z_k$ at times $0 < t_1 < \dots < t_k \dots \in [0,T]$
one would like to infer what they are.  But we can not observe
the entity $q$ directly---like in Plato's cave allegory we can only see
a `shadow' of it, formally given by a `measurement operator'
\begin{equation}  \label{eq:iI}
Y: \C{Q} \times \C{U} \ni (q,u(t_k)) \mapsto y_k = Y(q; u(t_k)) \in \C{Y};
\end{equation}
at least this is our model of what we are measuring.  
Frequently the space $\C{Y}$ may be regarded as finite dimensional,
as one can anly observe a finite number of quantities.
Usually the observation
will deviate from what we expect to observe even if we knew the right $q$
as \refeq{eq:I} is only a \emph{model}---so there is some model error
$\epsilon$, and the measurement will be polluted by some measurement error
$\vepsilon$.  Hence we observe $z_k = y_k + \epsilon + \vepsilon$.
From this one would like to know what $q$ and $u(t_k)$ are.  For the
sake of simplicity we will only consider one error term
$z_k = y_k + \vepsilon$ which subsumes all the errors.

The mapping in \refeq{eq:iI} is usually not invertible and hence the problem
is called ill-posed.  One way to address this is via regularisation
(see e.g.\ \cite{Engl2000}),
but here we follow a different track.  Modelling our lack-of-knowledge
about $q$ and $u(t_k)$ in a Bayesian way \cite{Tarantola2004} by replacing them
with a $\C{Q}$- resp.\ $\C{U}$-valued random variable (RV), the problem
becomes well-posed \cite{Stuart2010}.  But of course one is looking now at the
problem of finding a probability distribution that best fits the data;
and one also obtains a probability distribution, not just \emph{one} pair
$q$ and $u(t_k)$.  Here we focus on the use of a linear Bayesian approach 
\cite{Goldstein2007} in the framework of `white noise' analysis. 

The mathematical setup then is as follows: we assume that $\vOmega$
is a measure space with $\sigma$-algebra $\F{A}$ and
with a probability measure $\D{P}$, and that
$q: \vOmega \to \C{Q}$ and $u: \vOmega \to \C{U}$ are random variables.
For simplicity, we shall also require $\C{Q}$ to be a Hilbert space
where each vector is a possible realisation.  This is in order to allow
to measure the distance between different $q$'s as the norm of their
difference, and to allow the operations of linear algebra to be
performed.

\subsection{Bayesian updating of the measure} \label{SS:bayes-up-meas}
Bayes's theorem is commonly accepted as a consistent way to incorporate
new knowledge into a probabilistic description \cite{jaynes03,Tarantola2004}.
The elementary textbook statement of the theorem is about
conditional probabilities
\begin{equation}  \label{eq:iII}
 \D{P}(I_q|M_y) = \frac{\D{P}(M_y|I_q)}{\D{P}(M_y)}\D{P}(I_q),
\end{equation}
where $I_q$ is some subset of possible $q$'s, and $M_y$ is the information
provided by the measurement.  This becomes problematic when
the set $M_y$ has vanishing probability measure, but if all measures
involved have probability density functions (pdf), it may be formulated
as (\cite{Tarantola2004} Ch.\ 1.5)
\begin{equation}  \label{eq:iIIa}
 \pi_q(q|y) = \frac{p(y|q)}{\varpi} p_q(q),
\end{equation}
where $p_q$ is the pdf of $q$, $p(y|q)$ is the likelihood of $y=Y(q)$
given $q$, as a function of $q$ sometimes denoted by $L(q)$, and $\varpi$
is a normalising factor such that the conditional density $\pi_q(\cdot|y)$
integrates to unity.  These terms are in direct correspondence with
those in \refeq{eq:iII}.  Most computational approaches determine
the pdfs \cite{Marzouk2007,Stuart2010,Kucherova10}.

However,  Kolmogorov already
defined conditional probabilities via conditional expectation, e.g.\ see 
\cite{Bobrowski2006/087}.  Given the conditional expectation 
$\EXP{\cdot|M_y}$, the conditional
probability is easily recovered as $\D{P}(I_q|M_y) = \EXP{\chi_{I_q}|M_y}$,
where $\chi_{I_q}$ is the characteristic function of the subset $I_q$.

\subsection{Conditional expectation} \label{SS:cond-exp}
The easiest point of departure for conditional expectation in our setting
is to define it not just for one piece of measurement $M_y$, but for 
sub-$\sigma$-algebras $\F{S} \subset \F{A}$.  
The connection with an event $M_y$ is
to take $\F{S}=\sigma(Y)$, the $\sigma$-algebra generated by $Y$.

For RVs with finite variance---elements of $\C{S}:=L_2(\vOmega, \F{A},
\D{P})$---the conditional expectation $\EXP{\cdot | \F{S}}$ is defined as
the orthogonal projection onto $L_2(\vOmega, \F{S}, \D{P})$.  It
can then be extended as a contraction onto all $L_p(\vOmega, \F{A}, \D{P})$
with $p\geq 1$, e.g.\ see \cite{Bobrowski2006/087}.

Let us define the space $\E{Q} := \C{Q} \otimes \C{S}$ of $\C{Q}$-valued RVs
of finite variance, and set $\E{Q}_n := \C{Q} \otimes \C{S}_n$ with
$\C{S}_n := L_2(\vOmega, \F{S}, \D{P})$ for the $\C{Q}$-valued RVs with finite
variance on the sub-$\sigma$-algebra $\F{S}$, representing the new information.

The Bayesian update as conditional expectation is now simply formulated: 
\begin{equation} \label{eq:iIII}
  \EXP{q |\sigma(Y)} = P_{\E{Q}_n}(q) = \text{arg min}_{\tilde{q}\in\E{Q}_n}
  \| q - \tilde{q} \|^2_{\E{Q}},
\end{equation}
where $P_{\E{Q}_n}$ is the orthogonal projector onto $\E{Q}_n$.
Already in \cite{Kalman} it was noted that the conditional expectation
is the best estimate not only for the \emph{loss function} `distance squared',
as in \refeq{eq:iIII},
but for a much larger class of loss functions under certain distributional
constraints.  However for the above loss function this is valid without
any restrictions.

Requiring the derivative of the loss function in \refeq{eq:iIII} to
vanish---equivalently remembering from elementary geometry that the
line to the closest point is perpendicular to the approximating
subspace---one arrives at the Galerkin orthogonality conditions
\begin{equation}  \label{eq:iIV}
  \forall \tilde{q} \in \E{Q}_n: \quad
  \ip{q  - \EXP{q|\sigma(Y)}}{\tilde{q}}_{\E{Q}} = 0.
\end{equation}

To continue, note that the \emph{Doob-Dynkin} lemma
\cite{Bobrowski2006/087} assures us that if a RV like $\EXP{q|\sigma(Y)}$ is
measurable w.r.t.\ $\sigma(Y)$, then $\EXP{q|\sigma(Y)} = \psi(Y)$
for some measurable $\psi\in L_0(\C{Y};\C{Q})$.  More precisely one should
write $\EXP{q|\sigma(Y)} = \psi(Y(q)) = \psi \circ Y \circ q$.

Hence
$\E{Q}_n =  \overline{\text{span}} \{\phi \circ Y \circ q \in \E{Q} 
 \; | \; \phi \in L_0(\C{Y};\C{Q})\}$, where
$L_0(\C{Y};\C{Q})$ is the vector space of measurable maps from $\C{Y}$
to $\C{Q}$.  In particular one sees that $\EXP{q|\sigma(Y)}$ is of
this form.  In this light the task of finding
the conditional expectation may be seen as rephrasing
\refeq{eq:iIII} as: find  $\psi\in L_0(\C{Y};\C{Q})$ such that
\begin{equation} \label{eq:iV}
  \psi = \text{arg min}_{\phi\in L_0(\C{Y};\C{Q})}
  \|q - \phi \circ Y \circ q \|^2_{\E{Q}}.
\end{equation}
Then $q_a := \psi(y)=\EXP{q|y}$ is called the \emph{analysis},
\emph{assimilated}, or \emph{posterior} value,
incorporating the new information.

We would like to emphasise that it is the vector space setting of
$\C{Q}$ and $\C{Y}$ which has made this well-known formulation possible
\cite{Kalman},
and it will also allow for easy numerical computation.  To work with
measures as in \refeq{eq:iII} is cumbersome, as probability measures
are on the intersection of the unit sphere and the positive cone in
the space of signed finite measures.  
Similarly, in \refeq{eq:iIIa} the pdfs are in the positive cone of
$L_1(\vOmega)$ and on the unit sphere in $L_1(\vOmega)$, as they
have to integrate to unity.
A bit easier would be to work with
RVs which are in a metric space, the so-called Fr\'echet-type 
conditional expectation then minimises
the metric distance squared; but the Hilbert space setting is certainly
the simplest instance of this, and we will adhere to it here.

In case the $q$'s are not without constraints, or not in a vector space,
then they should be mapped to such quantities.
For example, if $q$ is a diffusion
tensor field, then it has to be symmetric and positive definite.
The symmetric tensors are of course a subspace, but the sub-manifold of
positive definite ones is not a subspace, nor is this subset closed,
making a minimisation as in \refeq{eq:iIII} ill-defined.
However the symmetric positive definite tensors can be
given the structure of a Lie group and a Riemannian manifold \cite{Arsigny_SPDM},
and then distance is measured as the length of a path along a geodesic.
Furthermore, the associated Lie algebra---the tangent space at the neutral element
of the Lie group---is in one to one correspondence with the geodesics; hence
one can play everything back to a vector space.  A simple case of this are
positive scalars; through the logarithm they are transformed into
a vector space without constraints.

\subsection{The linear Bayesian update} \label{SS:bayes-lin}
As we work in a vector space, we make
an approximation to simplify the computations by replacing
$L_0(\C{Y};\C{Q})$ above by $\E{L}(\C{Y},\C{Q}) \subset L_0(\C{Y};\C{Q})$,
the subspace of linear
continuous maps.  The minimisation \refeq{eq:iV} is then translated to:
find $K \in \E{L}(\C{Y},\C{Q})$ such that \cite{Kalman,Luenberger1969}
\begin{equation} \label{eq:iVI}
  K = \text{arg min}_{H\in \E{L}(\C{Y},\C{Q})}
    \| q - H \circ Y \circ q \|^2_{\E{Q}},
\end{equation}
and we set $\EXP{q|\sigma(Y)}_\ell := K \circ Y \circ q$, a linear in $Y$
approximation to $\EXP{q|\sigma(Y)}$.
As the projection is now onto the smaller space $\E{Q}_\ell :=
\overline{\text{span}} \{H \circ Y \circ q \in \E{Q}\; | \;  
H \in \E{L}(\C{Y},\C{Q})\} \subset \E{Q}_n \subset \E{Q}$, 
we are not using all the
information available.  Hence the error---the value of the functional
being minimised---will remain larger, but the computation is simpler.
The optimal $K$ is not hard to find by taking the derivative in \refeq{eq:iVI}
w.r.t.\ the linear map $H$ (see  e.g.\ \cite{Kalman} and \cite{Luenberger1969}
Ch.\ 3.2 Thm.\ 4.711), 
equivalently one may rewrite the variational Galerkin orthogonality
condition \refeq{eq:iIV}, for the final result see \refeq{eq:iVII}.

In the case of prior information represented by the \emph{forecast} RV
$q_f$, which results in the measurement forecast $y_f = Y(q_f)$,
the projection is adjusted by an affine shift to \cite{Kalman,Luenberger1969}
\begin{equation}  \label{eq:iVII}
  q_a = q_f + K(z - y_f), \text{   with   } 
  K = C_{q,y} (C_y + C_\vepsilon)^{-1},
\end{equation}
here the operator $K$ is also known as the \emph{Kalman gain}.
This includes the errors $\vepsilon$ assumed independent of $q$ with zero
mean and covariance operator 
$C_\vepsilon$, where $C_y := \EXP{\tilde{Y}(q,u)\otimes \tilde{Y}(q,u)}$ and
$C_{q,y} = \EXP{\tilde{q} \otimes \tilde{Y}(q,u)}$, where for any RV like $q$
for the sake of
brevity we set $\bar{q} := \EXP{q}$ such that $\tilde{q} := q - \bar{q}$
is the zero-mean part.
In case $C_y + C_\vepsilon$ is not invertible or close to singularity,
its inverse in \refeq{eq:iVII} should be replaced by the Moore-Penrose
pseudo-inverse.  This update is in some ways very similar to the
`Bayes linear' approach \cite{Goldstein2007}.  If the mean is taken in
\refeq{eq:iVII}, one obtains the familiar Kalman filter formula \cite{Kalman}
for the update of the mean, and one may show \cite{opBvrAlHgm11}
that \refeq{eq:iVII} also contains the Kalman update for the covariance.

Stated differently, in the situation of \refeq{eq:iVII} $q_a$ is the 
orthogonal projection of $q$ onto the subspace $\E{Q}_a = \E{Q}_f + \E{Q}_\ell$,
which is generated jointly by the prior information and the measurement.
This may be written as $\E{Q}_a = \E{Q}_f + \E{Q}_\ell = \E{Q}_f \oplus
\E{Q}_i$, where the information gain or innovation $\E{Q}_i$ is the part
of $\E{Q}_\ell$ orthogonal to $\E{Q}_f$, the last orthogonal sum reflecting
the terms in \refeq{eq:iVII}.

Each new measurement enlarges the space we project onto, and thus constrains
the error $ q − q_a$---which is in the orthogonal complement---further.  
In case the space generated by the measurements is not dense in $\E{Q}$
a residual error will thus remain, as the measurements do not contain enough
information to resolve our lack of knowledge about $q$.  Anyway,
finding $q$ is limited by the presence of the error $\vepsilon$, 
as obviously the error influences the update in \refeq{eq:iVII}.
If the measurement operator is approximated in some way---as it will be
in the computational examples to follow---this will introduce a new error,
further limiting the resolution.

%
%
%
%


%

\section{Numerical realisation} \label{S:num-real}
In the instances where we want to employ the theory detailed in the
previous \refS{bayes}, the spaces $\C{U}$ and $\C{Q}$ are usually infinite
dimensional, as is the space $\C{S} = L_2(\vOmega)$.
For an actual computation they have to be discretised or approximated
by finite dimensional spaces.  In our examples we will chose finite
element discretisations and corresponding subspaces.  Hence let
$\C{Q}_M := \text{span }\{\vrho_m\ : m=1,\dots,M\} \subset \C{Q}$ be an
$M$-dimensional subspace with basis $\{\vrho_m\}_{m=1}^M$.  An element
of $\C{Q}_M$ will be represented by the vector $\vek{q}=[q^1, 
\dots, q^M]^T \in \D{R}^M$ such that $\sum^M_{m=1} 
q^m \vrho_m \in \C{Q}_M$.  The space of possible measurements
can usually be taken to be finite dimensional, whose elements similarly
are represented by a vector of coefficients $\vek{z} \in \D{R}^R$.

On $\D{R}^M$, representing $\C{Q}_M$, the minimisation in \refeq{eq:iVI}
is translated to
\begin{equation} \label{eq:iVI-M}
  \vek{K} = \text{arg min}_{\vek{H} \in \D{R}^{M \times R}}
    \nd{\vek{q} - \vek{H}(\vek{Y}(\vek{q}))}^2_{\vek{Q}},
\end{equation}
where the mapping induced by $Y$ has been denoted by boldface $\vek{Y}$,
and the norm $\nd{\vek{q}}_{\vek{Q}}$ results from the inner product
$\bkt{\vek{q_1}}{\vek{q_2}}_{\vek{Q}} := \EXP{\vek{q_1}^T \vek{Q}\vek{q_2}}$
with $\vek{Q}_{m n} = \bkt{q_m}{q_n}_{\E{Q}}$, the Gram matrix
of the basis.  We  will later choose an orthonormal basis, so that
$\vek{Q} = \vek{I}$ is the identity matrix.  Then the update corresponding
to \refeq{eq:iVII} is
\begin{equation}  \label{eq:iIX}
  \vek{q}_a = \vek{q}_f + \vek{K}(\vek{z} - \vek{y}_f), \text{   with   } 
  \vek{K} = \vek{C}_{q,y} (\vek{C}_y + \vek{C}_\vepsilon)^{-1}.
\end{equation}
Here the covariances are naturally $\vek{C}_{q,y} := \EXP{\tilde{\vek{q}}\;
 \tilde{\vek{y}}^T} = \EXP{\tilde{\vek{q}} \otimes \tilde{\vek{y}}}$,
and similarly for $\vek{C}_y$ and $\vek{C}_\vepsilon$.

\subsection{Markov chain Monte Carlo} \label{SS:mcmc}
We shall shortly sketch the Markov chain Monte Carlo (MCMC) method for the
sake of completeness \cite{Madras-Fields:2002,Gamerman06}, as it will
be used on some examples.  It is a method which changes the underlying
probability measure according to \refeq{eq:iII} and \refeq{eq:iIIa} in
\refSS{bayes-up-meas}.

The idea is to construct a Markov chain with an equilibrium distribution
which corresponds to some desired distribution.  In our case this will
be the posterior or conditional distribution.  The method has the
distinct advantage that it does not need target probabilities but only
ratios of target probabilities to work.  This means that the normalisation
constant $\varpi$ in \refeq{eq:iIIa} does not have to computed at
all---this would involve a cumbersome high-dimensional integration---as
it cancels out.  Having such a Markov Chain, all that is required to
obtain samples from the posterior is to run the Markov chain for a
sufficiently long time.  As the posterior is the equilibrium distribution,
one is indeed sampling the posterior.

The Metropolis algorithm to achieve the construction is described very
quickly for the simplest case \cite{Madras-Fields:2002}:
assume that the datum $y$ has been observed, and that the range of possible
$q$'s has been quantised into $X$ equal bins with representatives
$\{q_\xi\}_{\xi=1}^X$, to each $q_\xi$ assign the number 
$\rho_\xi := L(q_\xi) p_q(q_\xi)$---the product of likelihood and 
prior pdf.  The representatives of the bins $q_\xi$ will be the states of 
the Markov chain, which will be denoted by $\{\vsigma_k\}_{k=1,\dots}$
in the order visited at step $k=1,\dots$.

For each step $k=1,2,\dots$ do the following: 

if currently in state  $\vsigma_k=q_\xi$, then 
\begin{enumerate}
\item pick any state $q_\zeta (\zeta \ne \xi)$ randomly
      with probability $1/(X-1)$,\\
      this is the \emph{proposal};
\item let $\alpha = \min\{1, \rho_\zeta / \rho_\xi\}$, this is the
      \emph{acceptance probability};
\item accept $q_\zeta$ with probability $\alpha$,\\ i.e.\ pick a sample
      of a uniformly distributed RV $U \in [0,1]$, and \\
      if $U \le \alpha$ then
      set $\vsigma_{k+1} =q_\zeta$,\\
      else set $\vsigma_{k+1}=q_\xi$.
\end{enumerate}
As the acceptance probability in step 2 only involves ratios of the $\rho_\xi$,
the---hard to compute---normalising constant $\varpi$ in \refeq{eq:iIIa}
is indeed not needed.

The simplicity of the formulation and possibility to correctly estimate
the posterior without any approximation are certainly the main advantages
of the method.
We should warn that although the formulation is so simple,
to use the method efficiently and correctly may not be so simple.  

First, MCMC is a Monte Carlo method, which means that estimates converge
only slowly with increasing number of samples.
Second, the samples---coming from a Markov chain---are obviously not 
independent.  This makes estimating any statistic other than the mean
(e.g.\ the variance) difficult.  Also usual statistical formulas
for the spread or accuracy of the estimates assume \emph{independent}
samples and are not applicable.
Third, the sequence of visited states is not even \emph{stationary} as
the equilibrium distribution is only the asymptotic distribution.
This means that one has to wait until $k$ is \emph{large enough}
before actually starting to sample, this is called the \emph{burn-in}
period.  Usually estimating that the burn-in period is over has to be
tested by tests on the stationarity of the sequence.  
And finally, if the likelihood function and prior pdf differ very much,
the acceptance probabilities $\alpha$ in step 2 may be very low---meaning that
practically the chain has got stuck in some state and does not move on.

Here we have used one of the simplest variants, where the underlying
transition probabilities of the Markov chain are all equal.  For extensions
which partly try to alleviate the problems alluded to above
and thorough explanations the reader should refer to the literature cited
and the references therein.

Now, the two methods to be considered in the following are based on the
update via conditional
expectation as outlined in \refSS{cond-exp}, and they both employ
the linear approximation from \refSS{bayes-lin}.  This means that in contrast
to MCMC they make an additional approximation error, or stated differently,
they do not use all the information availabe for the update.
After the discretisation
of the space $\C{Q}$ this boils down to using \refeq{eq:iIX} in some way
on $\C{Q}_M$.  How this is done differs in the two methods.

\subsection{Ensemble Kalman Filter} \label{SS:enkf}
The ensemble Kalman filter (EnKF) is a Monte Carlo method, a sampling
interpretation of \refeq{eq:iIX}, for details see
e.g.\ \cite{Evensen2009,Evensen2009a}.
The basic form is actually simple to describe and is a good preparation
for the PCE-based form in the next \refSS{pce-bayes}.

Pick $Z$ independent samples $\{\vek{q}_f(\omega_z)\}_{z=1}^Z$ according 
to the probability measure $\D{P}$.  Now \refeq{eq:iIX} holds for the
the RV $\vek{q}: \vOmega \to \C{Q}_M$, i.e.\ for all $\omega \in \vOmega$
(a.s.---almost $\D{P}$-surely).  One takes this to mean that it holds for
each $\omega_z \in \vOmega$.  Arranging the the samples in a matrix
$\opb{Q}_f := [\vek{q}_f(\omega_1), \dots, \vek{q}_f(\omega_Z)]$, and
similarly for the forecasts $\opb{Y}_f := [\vek{y}_f(\omega_1),\dots,
\vek{y}_f(\omega_Z)]$ and measurements $\opb{Z}$, \refeq{eq:iIX} is now
formulated in matrix notation as
\begin{equation}  \label{eq:enkf-1}
  \opb{Q}_a = \opb{Q}_f + \vek{K}(\opb{Z} - \opb{Y}_f),
\end{equation}
where $\vek{K}$ is the same as in \refeq{eq:iIX}.

But the variances needed to compute $\vek{K}$ have to be estimated from the
sample.  This simply takes the form
\begin{eqnarray}  \label{eq:enkf-2}
  \vek{C}_{q,y} &\approx& \frac{1}{Z-1} \tilde{\opb{Q}}_f \tilde{\opb{Y}}_f^T\\
  \vek{C}_{y} &\approx& \frac{1}{Z-1} \tilde{\opb{Y}}_f \tilde{\opb{Y}}_f^T.
  \label{eq:enkf-3}
\end{eqnarray}
The normalisation terms $(Z-1)^{-1}$ do not really have to be used as they
cancel in the computation of $\vek{K}$.  With the estimated means 
$$ \bar{\vek{q}}_f = \frac{1}{Z} \sum_{z=1}^Z \vek{q}_f(\omega_z) \qquad \text{and}
\qquad \bar{\vek{y}}_f = \frac{1}{Z} \sum_{z=1}^Z \vek{y}_f(\omega_z),$$
the terms in \refeq{eq:enkf-2} and \refeq{eq:enkf-3} are
\begin{eqnarray}  \label{eq:enkf-4}
  \tilde{\opb{Q}}_f &=& \opb{Q}_f - \bar{\vek{q}}_f\, \vek{1}^T_Z\\
  \tilde{\opb{Y}}_f &=& \opb{Y}_f - \bar{\vek{y}}_f\, \vek{1}^T_Z,
  \label{eq:enkf-5}
\end{eqnarray}
where $\vek{1}_Z$ is a vector of ones of size $Z$.

This method is a Monte Carlo method, hence it also suffers from the slow
convergence with increasing $Z$.  On the other hand it is fairly simple
to implement, all it needs are samples.  In practice the number of
samples is often low, and then special care is needed when computing the
covariances---whose rank can not be higher than the number of samples---and
the Kalman gain $\vek{K}$, see \cite{Evensen2009,Evensen2009a}.
Another closely connected problem is the possible underestimation of
the posterior sample variance with low sample numbers.
Still, this method is currently
a favourite for problems where the computation of the predicted measurement
$\vek{y}_f(\omega_z)$ is difficult or expensive.  It needs far fewer samples
for meaningful results than MCMC, but on the other hand it uses the linear
approximation inherent in \refeq{eq:iIX}.

\subsection{Polynomial chaos projection}  \label{SS:pce-bayes}
In \refSS{enkf} the \refeq{eq:iIX} was discretised in the variables
$\omega \in \vOmega$ through sampling.  Here another track is taken,
relying on the tensor product representation of the RVs
$\E{Q} = \C{Q} \otimes \C{S}$ from \refSS{cond-exp}.  The first
factor already was discretised to $\C{Q}_N \subset \C{Q}$, and
here an explicit discretisation of the second factor is used.
In a numerical sense the sampling in \refSS{enkf} is of course also
a discretisation.  As for $\C{Q}_M$ we pick a finite set of independent
vectors in $\C{S}$.  As $\C{S} = L_2(\vOmega)$, these abstract vectors
are in fact RVs with finite variance.  Here we will use \emph{Wiener}'s 
\emph{polynomial chaos} expansion (PCE) as basis 
\cite{Janson1997,matthies6}, this allows
to use \refeq{eq:iIX} without sampling, see 
\cite{opBvrAlHgm11,bvrAlOpHgm11,boulder:2011}, and also
\cite{saadGhn:2009,Blanchard2010a}.

The PCE is an expansion in multivariate \emph{Hermite polynomials} 
\cite{Janson1997,matthies6}; denote
by $H_{\vek{\alpha}}(\vek{\theta}) = \prod_{k \in \D{N}} h_{\alpha_k}(\theta_k)
\in \C{S}$ the multivariate polynomial in standard independent Gaussian RVs 
$\vek{\theta}(\omega) = (\theta_1(\omega),\dots, \theta_k(\omega), 
\dots)_{k\in \D{N}}$,
where $h_j$ is the usual univariate Hermite polynomial, and $\vek{\alpha} = 
(\alpha_1,\dots,\alpha_k,\dots)_{k\in \D{N}}\in \C{N}:=\D{N}_0^{(\D{N})}$
is a multi-index of generally infinite lenght but with only finitely many
entries non-zero.  As $h_0 \equiv 1$,
the infinite product is effectively finite and always well-defined.

The \emph{Cameron-Martin} theorem assures us \cite{Janson1997}
that the set of these polynomials is dense in $\C{S} = L_2(\vOmega)$,
and in fact $\{H_{\vek{\alpha}}/\sqrt{(\vek{\alpha} !)} \}_{\vek{\alpha} 
\in \C{N}}$ is a complete orthonormal system (CONS), 
where $\vek{\alpha} ! := \prod_{k \in \D{N}} (\alpha_k !)$ is the product
of the individual factorials, also well-defined as except for finitely many 
$k$ one has $\alpha_k ! = 0! = 1$.  So we may assume that $\vek{q}(\omega) =
\sum_{\vek{\alpha}\in \C{N}} \vek{q}^{\vek{\alpha}} H_{\vek{\alpha}}
(\vek{\theta}(\omega))$, and similarly for $\vek{z}$ and $\vek{y}$.  In this
way the RVs are expressed as functions of other, known RVs, and not
through samples.

The space $\C{S}$ may now be discretised by taking a finite subset $\C{J}
\subset \C{N}$ of size $J = \ns{\C{J}}$, and setting $\C{S}_J = \text{span }
\{H_{\vek{\alpha}}\,:\, \vek{\alpha} \in \C{J} \} \subset \C{S}$.  The
orthogonal projection $P_J$ onto $\C{S}_J$ is then simply
\begin{equation}  \label{eq:proj-J}
P_J: \C{Q}_M \otimes \C{S} \ni
\sum_{\vek{\alpha}\in \C{N}} \vek{q}^{\vek{\alpha}} H_{\vek{\alpha}} \mapsto
\sum_{\vek{\alpha}\in \C{J}} \vek{q}^{\vek{\alpha}} H_{\vek{\alpha}} 
\in \C{Q}_M \otimes \C{S}_J.
\end{equation}
We then take \refeq{eq:iIX} and rewrite it as
\begin{eqnarray}  \label{eq:proj-lin-f1}
  \vek{q}_a &=& \vek{q}_f + \vek{K}(\vek{z} - \vek{y}_f) =\\
  \sum_{\vek{\alpha}\in \C{N}} \vek{q}_a^{\vek{\alpha}} H_{\vek{\alpha}} &=& 
  \sum_{\vek{\alpha}\in \C{N}} \left(\vek{q}_f^{\vek{\alpha}} + \vek{K}\left( 
  \vek{z}^{\vek{\alpha}}-\vek{y}^{\vek{\alpha}}_f\right)\right)H_{\vek{\alpha}}.
  \label{eq:proj-lin-f2}
\end{eqnarray}
Projecting both sides of \refeq{eq:proj-lin-f2} is then very
simple and results in
\begin{equation} \label{eq:proj-lin-J}
  \sum_{\vek{\alpha}\in \C{J}} \vek{q}_a^{\vek{\alpha}} H_{\vek{\alpha}} = 
  \sum_{\vek{\alpha}\in \C{J}} \left(\vek{q}_f^{\vek{\alpha}} + \vek{K}\left( 
  \vek{z}^{\vek{\alpha}}-\vek{y}^{\vek{\alpha}}_f\right)\right)H_{\vek{\alpha}}.
\end{equation}
Obviously the projection $P_J$ commutes with the Kalman operator $K$ and
hence with its finite dimensional analogue $\vek{K}$.  One may actually
concisely write \refeq{eq:proj-lin-J} as
\begin{equation} \label{eq:proj-comm-K}
  P_J \vek{q}_a = P_J \vek{q}_f + P_J \vek{K}(\vek{z} - \vek{y}_f) =
  P_J\vek{q}_f + \vek{K}(P_J\vek{z} - P_J\vek{y}_f).
\end{equation}

Elements of the discretised space $\E{Q}_{M,J} = \C{Q}_M \otimes
\C{S}_J \subset \E{Q}$ thus may be written as 
$\sum_{m=1}^M \sum_{\vek{\alpha}\in \C{J}} q^{\vek{\alpha},m} \vrho_m 
H_{\vek{\alpha}}$.  The tensor representation 
$\mat{q} := (q^{\vek{\alpha},m}) = \sum_{\vek{\alpha}\in \C{J}} 
\vek{q}^{\vek{\alpha}} \otimes \vek{e}^{\vek{\alpha}}$, where the 
$\vek{e}^{\vek{\alpha}}$ are the unit vectors in $\D{R}^J$, may be used
to express \refeq{eq:proj-lin-J} or \refeq{eq:proj-comm-K} succinctly as
\begin{equation}  \label{eq:proj-t}
 \mat{q}_a = \mat{q}_f+ \mat{K}(\mat{z}-\mat{y}_f),
\end{equation}
where $\mat{K} = \vek{K} \otimes \vek{I}$ with $\vek{K}$ from \refeq{eq:iIX}.
Hence the update equation is naturally in a tensorised form.  This is
how the update can finally be computed in the PCE representation without
any sampling \cite{opBvrAlHgm11,bvrAlOpHgm11,boulder:2011}.

It remains to say how to compute the Kalman gain---or rather the covariance
matrices---in this approach.  Given the PCEs of the RVs, this is actually
quite simple as any moment can be computed directly from the PCE
\cite{matthies6,opBvrAlHgm11,bvrAlOpHgm11}; remembering that 
$\vek{q}^0=\EXP{\vek{q}}$, due to the orthogonality of the 
$H_{\vek{\alpha}}$ one has for the variance
\begin{eqnarray}  \label{eq:cov-PCE-1}
   \vek{C}_{y} 
   &\approx& \vek{C}_{P_J y} = \EXP{(P_J\vek{y}) \otimes (P_J\vek{y})} =
   \sum_{\vek{\alpha}\in \C{J}, \vek{\alpha} \ne 0} (\vek{\alpha} !)\;
   \vek{y}^{\vek{\alpha}}\otimes \vek{y}^{\vek{\alpha}},\\
  \vek{C}_{q,y} &=& 
  \sum_{\vek{\alpha}\in \C{N}, \vek{\alpha} \ne 0} (\vek{\alpha} !)\;
   \vek{q}^{\vek{\alpha}}\otimes \vek{y}^{\vek{\alpha}} \approx
  \sum_{\vek{\alpha}\in \C{J}, \vek{\alpha} \ne 0} (\vek{\alpha} !)\;
   \vek{q}^{\vek{\alpha}}\otimes \vek{y}^{\vek{\alpha}}. 
    \label{eq:cov-PCE-2}
\end{eqnarray}

As was already remarked
\refSS{enkf}, the covariance matrix can not have higher rank than
the number of tensor products in the sum.  And certainly any such truncation
as is implicit in the projection operator $P_J$ \refeq{eq:proj-J} or
explicit as in \refeq{eq:cov-PCE-1} and \refeq{eq:cov-PCE-2} will reduce
the variance.  The difference to the procedure in \refSS{enkf} is that
here the $\vek{q}^{\vek{\alpha}}$ and $\vek{y}^{\vek{\alpha}}$ are
components of the orthogonal basis $H_{\vek{\alpha}}$.
The main `difficulty' to apply the the update as described in this section
is the need to have a PCE of $\vek{q}_f(\omega)$ and especially of the
forecast measurement $\vek{y}_f(\omega) = Y(\vek{q}_f(\omega))$.  How
this might be done will be sketched in the next \refS{forward}.

One should not forget to draw attention to the fact that once such
an expansion is available, evaluating an expression such as
$\sum_{\vek{\alpha}\in \C{J}} \vek{y}^{\vek{\alpha}}
          H_{\vek{\alpha}}(\omega)$
for some specific $\omega \in \vOmega$ is computationally relatively cheap,
and this is an approximation to a sample of $\vek{y}_f(\omega)$---it is
precisely $P_J \vek{y}_f(\omega)$.
Hence this can be used wherever samples are needed, and this device has
actually been employed, i.e.\ in the MCMC method 
\cite{Marzouk2009a,Kucherova10,Pence2010,kucerovaSykBvrHgm:2011}
as described in \refSS{mcmc}, or in the EnKF method \cite{Li2009}
described in \refSS{enkf}.

%
%
%
%
%
%


%

\section{The forward problem} \label{S:forward}
As was noted already in \refS{num-real}, the spaces $\C{U}$ and $\C{Q}$
are usually infinite dimensional, as is the space $\C{S} = L_2(\vOmega)$.
Similarly to the discretisation described there, the space $\C{U}$ has
to be discretised as well.  In an analogous fashion, choosing an
$N$-dimensional subspace $\C{U}_N = \text{span }\{\upsilon_n\ : n=1,\dots,N\}
\subset \C{U}$ with basis $\{\upsilon_n\}_{n=1}^N$.  An element
of $\C{U}_N$ will similarly be represented by the vector $\vek{u}=[u^1, 
\dots, u^N]^T \in \D{R}^N$ such that $\sum^N_{n=1} 
u^n \upsilon_n \in \C{U}_N$.

The state equation may then be discretised by inserting the `ansatz' that
the solution is from $\C{U}_N$ by assuming in \refeq{eq:I} that
$u = \sum_n u^n \upsilon_n$.  We obtain equations for the coefficients $u^n$
for example by projecting in general Galerkin manner \refeq{eq:I} onto
$\C{U}_N$; most straightforward is to use the basis $\{\upsilon_n\}_{n=1}^N$
also for the projection:
\begin{equation}  \label{eq:I-proj}  \forall k: \quad
\bkt{\upsilon_k}{\frac{\dd}{\dd t}\sum_n u^n(t) \upsilon_n}_{\C{U}} + 
\bkt{\upsilon_k}{A(q;\sum_n u^n(t) \upsilon_n)}_{\C{U}} = 
\bkt{\upsilon_k}{f(q;t)}_{\C{U}},
\end{equation}
where we have now only included the dependence on the parameter $q \in \C{Q}$
which has to be identified.
\refeq{eq:I-proj} may in standard manner succinctly be written as an equation
on $\D{R}^N$:
\begin{equation}  \label{eq:Ia-proj} 
\frac{\dd}{\dd t} \vek{u}(t) + \vek{A}(\vek{q};\vek{u}) = \vek{f}(\vek{q};t).
\end{equation}
The measurement operator will, in combination with the developments in
\refS{num-real}, be denoted as (see \refeq{eq:iI})
\begin{equation}  \label{eq:meas-disc}
\vek{y}_f = \vek{Y}(\vek{u}(t_k),\vek{q}_{k-1}) =  \vek{Y}(\vek{q}_{k-1}).
\end{equation}

The unknown parameter is modelled as a RV---see \refS{bayes}, which
in the discrete setting of \refS{num-real} means that \refeq{eq:Ia-proj}
reads
\begin{equation}  \label{eq:I-RV} 
\frac{\dd}{\dd t} \vek{u}(\omega) + \vek{A}(\vek{q}(\omega);\vek{u}(\omega))
 = \vek{f}(\vek{q}(\omega)) \quad \D{P}\text{- a.s. in } \omega \in \vOmega.
\end{equation}
For MCMC (\refSS{mcmc}) and EnKF (\refSS{enkf}) or any other Monte Carlo
or collocation method, this equation has to be solved for each sample 
or collocation point $\omega_z$ and $\vek{q}(\omega_z)$,
to obtain $\vek{u}(\omega_z)$, and with that predict the measurement 
$\vek{y}_f(\omega_z) = \vek{Y}(\vek{u}(\omega_z),\vek{q}(\omega_z))$.
Obviously, this may be computationally quite costly.
As remarked at the end of \refSS{pce-bayes}, if one can compute a PCE of
$\vek{y}_f(\omega) = \sum_{\vek{\alpha}\in \C{J}} 
\vek{y}^{\vek{\alpha}}_f H_{\vek{\alpha}}(\vek{\theta}(\omega))$,
this may be used instead to great advantage, as it will be usually much
less work to evaluate instead of going via \refeq{eq:I-RV}.

For the linear Bayesian update in \refSS{pce-bayes} via PCE this comes natural,
but of course this means that \refeq{eq:I-RV} has to be solved so that this
is possible.  Just as the parameter was represented in
\refSS{pce-bayes} via its truncated PCE \refeq{eq:proj-J} 
$\vek{q}(\omega) = \sum_{\vek{\alpha}\in \C{J}} \vek{q}^{\vek{\alpha}} 
H_{\vek{\alpha}}(\vek{\theta}(\omega))$, or through the tensor of
coefficients $\mat{q} = \sum_{\vek{\alpha}\in \C{J}} \vek{q}^{\vek{\alpha}} 
\otimes \vek{e}^{\vek{\alpha}}$, the same is assumed by an `ansatz'
for the solution to \refeq{eq:I-RV} $\vek{u}(\omega) = 
\sum_{\vek{\alpha}\in \C{J}} \vek{u}^{\vek{\alpha}} 
H_{\vek{\alpha}}(\vek{\theta}(\omega))$, represented through $\mat{u} = 
\sum_{\vek{\alpha}\in \C{J}} \vek{u}^{\vek{\alpha}} 
\otimes \vek{e}^{\vek{\alpha}}$.  

Inserting this into \refeq{eq:I-RV}
and projecting with $P_J$ \refeq{eq:proj-J} onto $\C{S}_J$, one obtains the
equations for $\mat{u}(t) \in \D{R}^{N} \otimes \D{R}^{J}$ through
Galerkin conditions with the basis $\{H_{\vek{\alpha}}\}_{\vek{\alpha}\in \C{J}}$
analogous to \refeq{eq:I-proj}:
\begin{equation}  \label{eq:I-tens} 
\frac{\dd}{\dd t} \mat{u} + \mat{A}(\mat{q};\mat{u}) = \mat{f}(\mat{q}),
\end{equation}
with the obvious interpretation of the terms.  The same procedure applied to
\refeq{eq:meas-disc} results in $\mat{y}_f = 
\mat{Y}(\mat{u}(t_k);\mat{q}_{k-1})$.  Identifying $\mat{q}_f$ with 
$\mat{q}_{k-1}$ at step $k$ of the updating, all the terms needed to
use the update \refeq{eq:proj-t} are present.

Let us remark that the tensor product representation like the
one employed for the state $\mat{u}$ in \refeq{eq:I-tens} may be
extended \cite{boulder:2011} to all entities in \refeq{eq:I-tens}.
Low-rank approximations to those entities in tensor representation
then become possible in different guises, as model reduction
\cite{BoyLeBYMadE:2010}, as so called `generalised spectral decomposition'  
\cite{Nouy2009}, during the solution process progressively as
`proper generalised decomposition' \cite{NouyACM:2010}, or in the iteration as
approximate, perturbed, or compressed iteration \cite{whBKhEET:2008,Zander10},
and this may lead to considerable numerical savings.  As \refeq{eq:proj-t} is
in such a format, these savings carry over to the Bayesian update.

%
%
%
%
%
%


%

\section{Numerical examples} \label{S:num-xmpls}
A comparison of the different updating methods is outside the scope
of this paper and will be presented elsewhere.  Here we rather want
to present some examples of identification procedures which illustrate
how the methods work.  Note that these examples are not intended to
describe real world applications but
are rather set up to illustrate the class of full and 
linear Bayesian identification methods.  Specifically, the measurement
operation will not be preformed on a real system, but will also
be simulated through computation.  This has the advantage that we
know what the `truth' is, i.e.\ the value of the parameter field with
which we simulate the measurement, so that one can compare this with
the estimate from the identification.

These examples will be both linear and 
nonlinear problems described through partial differential equations,
namely linear and nonlinear diffusion and elasto-plastic behaviour.
Hence the \refeq{eq:I} is here a partial differential equation, and
the `parameter' to be identified will be some coefficient fields in
the equation, either the diffusion coefficient, or the shear modulus.
Let us denote such a generic random field as $q(x,\omega)$, where now
$x \in\E{G}$ is a point in the spatial domain $\E{G} \subset \D{R}^d$.  
More precisely the spatial function $q(\cdot,\omega)$ is the parameter to
be identified.  We assume that this field is transformed in such a way
that it is without constraints in some vector space $\C{Q}$.

To start the numerical computations, one more step is needed to
be able to effectively describe the random fields.  What one needs
for the developments described in \refSS{pce-bayes} is the PCE of
the random field $q(x,\omega) = \sum_{\vek{\alpha} \in \C{N}} 
q^{\vek{\alpha}}(x) H_{\vek{\alpha}}(\vek{\theta}(\omega))$.  The spatial
coefficient functions are given by simple projection $q^{\vek{\alpha}}(x)
= \EXP{q(x,\cdot) H_{\vek{\alpha}}(\vek{\theta}(\cdot))}$.  

The computational
problem with this approach is that the PCE is completely general and is
defined without any reference to the random field $q(x,\omega)$.
This means that an excessive number of RVs $\theta_1(\omega),
\dots,\theta_k(\omega),\dots$ may be needed to give an accurate enough
approximation when the above PCE is truncated to some $\vek{\alpha} \in 
\C{J} \subset \C{N}$.  

One way to have an accurate description with a
fairly small number of RVs is to use the \KL{} expansion (KLE), e.g.\ see
\cite{matthies6}.  With the covariance function of the random field
$c_q(x_1,x_2) := \EXP{q(x_1,\cdot) q(x_2,\cdot)}$---for the sake of
simplicity we will only consider scalar random fields---one computes
the eigenvalues and -functions of the Fredholm eigenvalue problem with
the covariance function as kernel:
\begin{equation}  \label{eq:fredh-KL}
  \int_{\E{G}} c_q(x_1,x_2) q_j(x_2) \; \di x_2 = \lambda_j \, q_j(x_1),
\end{equation}
which results in the spectral decomposition $c_q(x_1,x_2) = \sum_j 
\lambda_j q_j(x_1) q_j(x_2)$.

The \KL{} expansion then is
\begin{equation}   \label{eq:exp-KL}
  q(x,\omega) = \bar{q}(x) + \sum_j \sqrt{\lambda_j}\, \xi_j(\omega) q_j(x),
\end{equation}
where the non-negative eigenvalues $\lambda_j$ and orthonormal 
eigenfunctions $q_j(x)$ are usually ordered in a
decreasing sequence, and the zero-mean unit-variance uncorrelated
(i.e.\ orthonormal) RVs $\xi_j(\omega)$
are given through a simple projection of the random field
\begin{equation}   \label{eq:RV-KL}
\xi_j(\omega) = \frac{1}{\sqrt{\lambda_j}} 
      \int_{\E{G}} (q(x,\omega)-\bar{q}(x)) q_j(x)\; \di x.
\end{equation}
As the variance of each term in the KLE \refeq{eq:exp-KL}
is $\lambda_j$---the total variance is $\sigma^2_t := \sum_j \lambda_j = 
\int_{\E{G}} c_q(x,x)\, \di x$---it is usually truncated after
say $M$ terms so that the residual variance $\sigma^2_r := 
\sum_{j>M} \lambda_j$ is small.

The field $q(x,\omega)$ is then approximated or represented by
\begin{equation}   \label{eq:exp-KL-r}
  q(x,\omega) = \bar{q}(x) + \sigma_c \sum_{j \le M} 
    \sqrt{\lambda_j}\, \xi_j(\omega) q_j(x),
\end{equation}
which has the same mean $\bar{q}(x)$ as the original field \refeq{eq:exp-KL}
and covariance $c_q(x_1,x_2) = \sigma_c^2 \sum_{j \le M} \lambda_j 
q_j(x_1) q_j(x_2)$.  If the correction factor $\sigma_c$ is chosen as
$\sigma_c = \sqrt{\sigma_t^2/(\sigma_t^2 - \sigma_r^2)}$, the approximate
field \refeq{eq:exp-KL-r} will also have the correct total variance---only
distributed on the modes $q_j(x)$ with $j \le M$. 

Through the KLE the random field is best approximated in
variance with the least number of RVs.  For further computation the RVs
$\xi_j(\omega)$ are sometimes unwieldy, so they in turn may now be expanded
in a PCE
\begin{equation}  \label{eq:PCE-xi}
  \xi_j(\omega) = \sum_{\vek{\alpha} \in \D{N}_0^M} \xi^{\vek{\alpha}}_j
  H_{\vek{\alpha}}(\vek{\theta}_M(\omega)),
\end{equation}
where we only need $M$ isonormal RVs $\vek{\theta}_M(\omega) = 
(\theta_1(\omega),\dots,\theta_M(\omega))$.  To become computationally
viable, the series \refeq{eq:PCE-xi} will have to be reduced to a finite
sum, usually by choosing a finite subset $\C{J}_j \subset \D{N}_0^M$.
One criterion could of course be again to choose those $\vek{a}$ which
give the highest contribution.  By doing this one will obviously decrease
the variance, and this may be compensated by using an additional scaling
factor $\sigma_j$ computed similarly to $\sigma_c$ above.  For the sake
of simplicity we shall avoid doing this here.  Inserting the
thus truncated series into \refeq{eq:exp-KL-r} gives
\begin{equation}   \label{eq:exp-KL-red}
  q(x,\omega) = \bar{q}(x) + \sigma_c \sum_{j\le M}\sum_{\vek{\alpha}\in\C{J}_j}
    \sqrt{\lambda_j}\, \xi_j^{\vek{\alpha}} H_{\vek{\alpha}}
    (\vek{\theta}_M(\omega)) q_j(x).
\end{equation}
Let $\C{J} := \bigcup_{j} \C{J}_j$, set $\xi_j^{\vek{\alpha}} = 0$ for
$\vek{\alpha} \in \C{J} \setminus \bigcup_{k \ne j} \C{J}_k$, and let
$\xi^{\vek{\alpha}}(x) := \sigma_c \sum_{j\le M} \sqrt{\lambda_j} 
\xi_j^{\vek{\alpha}} q_j(x)$ and $\xi^0(x) =\bar{q}(x)$;
then the truncated PCE of $q(x,\omega)$ in $M$ Gaussian RVs is
\begin{equation}   \label{eq:exp-KL-PCE}
  q(x,\omega) = \sum_{\vek{\alpha}\in\C{J}}\, \xi^{\vek{\alpha}}(x)
   H_{\vek{\alpha}}(\vek{\theta}_M(\omega)).
\end{equation}
Depending on the specific situation, either the expression \refeq{eq:exp-KL-r}
or \refeq{eq:exp-KL-PCE} can be used in an actual computation.


\subsection{Examples with MCMC updating} \label{SS:diff}
Our study will start with a relatively simple system and then progress
to more difficult cases, illustrating the Markov chain Monte Carlo (MCMC)
methods of \refSS{mcmc}.  Let us recall that in that case the posterior
is described by changing the distribution of the underlying RVs.

\begin{figure}[h!]
\centering
\includegraphics[width=0.68\textwidth]{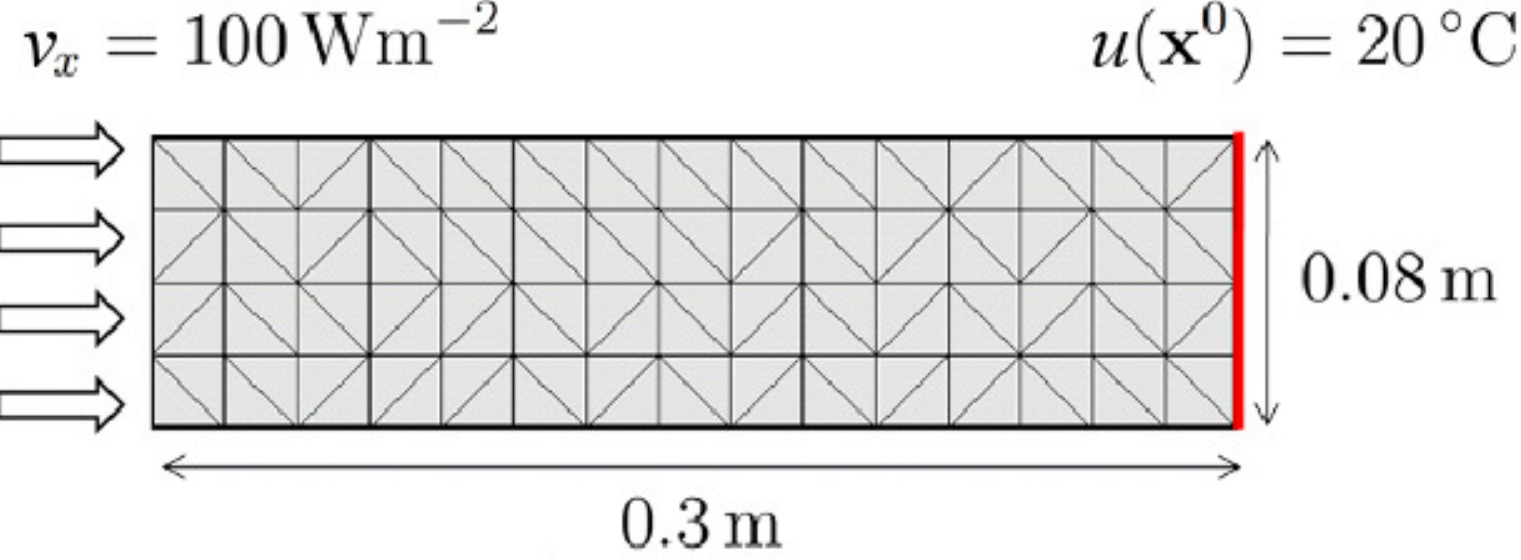}
\caption{Experimental setup.}
\label{F:schema}
\end{figure}

We start with a simple model of steady state heat transfer \cite{Kucherova10}, 
where the energy balance equation leads to
\begin{equation}   \label{eq:BAL1}
 -\divg (\kappa(x,\omega)\nabla u(x,\omega))=f(x,\omega), 
  \textrm{  a.e. } x\in \C{G} \subset \D{R}^2
\end{equation}
where both the heat conductivity $\kappa$ and the right hand side $f$ may be
considered as random fields over a probability space $\vOmega$.  
Thus \refeq{eq:BAL1} is required to hold almost surely in $\omega$, 
i.e.\ $\D{P}$-almost everywhere.  For the sake of simplicity the conductivity 
$\kappa(x,\omega)$  is taken to be a scalar valued random field, although a conductivity tensor would be a more appropriate choice.  First we will not let
$\kappa$ depend on $u$ and \refeq{eq:BAL1} represents linear model. 
To the balance equation one has to add appropriate Dirichlet and Neumann
boundary conditions, given as a a prescribed temperature
 $u_0=20^\circ \text{C}$ and a heat flux $\nu_x=100 \textrm{ Wm}^{-2}$
 on the left side and zero flux on top and bottom
respectively, as shown in \refig{schema}.
This completes the concrete description of the abstract \refeq{eq:I}, where
here there is no time dependence and thus the time derivative terms vanish.

After spatial discretisation by a standard finite element
method into $120$ elements (see \refig{schema}), the problem is further
discretised with the help of random variables as described in \refS{forward}. Following this the problem reduces to the form 
given by \refeq{eq:I-tens} without the first term as mentioned already.  
The system in \refeq{eq:I-tens} is then solved via the methods given
in \cite{BoyLeBYMadE:2010,NouyACM:2010,Zander10}.

The conductivity parameter $\kappa(x,\omega)$ is first considered in two 
different scenarios: first as a simple random variable---i.e.\ no dependence
on $x\in \C{G}$---and then as a random field.

\subsubsection{Conductivity as a simple random variable}\label{sec:condrv}
We take the thermal conductivity $\kappa$ as independent of the coordinates
$x\in \C{G}$, and for simplicity it is we assume it as an a priori normally
distributed variable $\kappa_f=q_f$ with mean value $\EXP{q_f} = 2
\text{ Wm}^{-1} \text{K}^{-1}$ and
standard deviation $0.3$ Wm$^{-1}$K$^{-1}$.

The measurements will be simulated with a `true' value of
$\kappa_t = 2.5\text{ Wm}^{-1} \text{K}^{-1}$.  They are 
performed in each node of the finite element mesh; and further polluted at
each node by independent centred Gaussian noise 
with a standard deviation of $\sigma_{\varepsilon} = 10^\circ\text{C}$. 

\begin{figure}[h!]
\centering
\includegraphics*[width=0.68\textwidth]{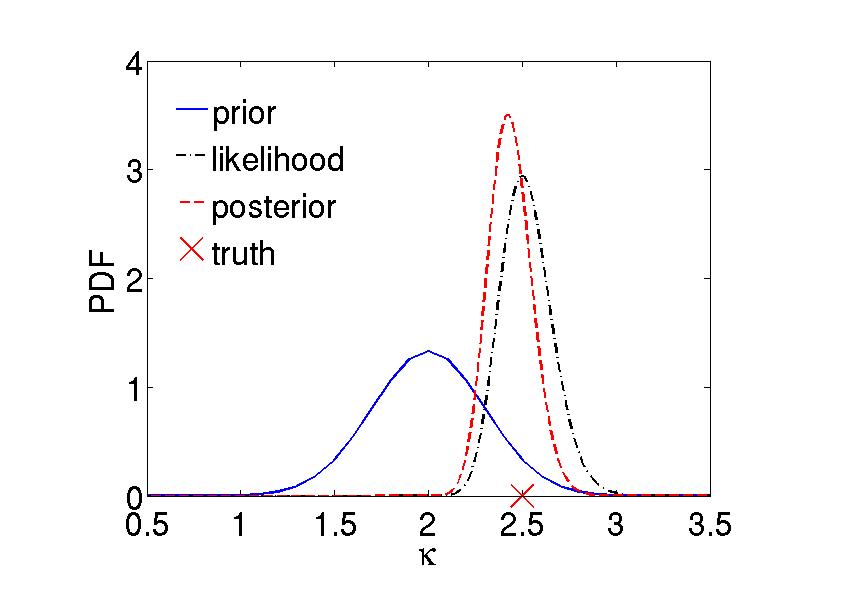}
\caption{Comparison of prior and posterior pdf and the likelihood function.}
\label{F:pdfcomp1}
\end{figure}
In this particular example the identification of $\kappa$ is done in a fully Bayesian manner (see \refeq{eq:iII}) with the help of a
MCMC procedure with $100000$ samples
as described in \refSS{mcmc}.  In \refig{pdfcomp1} we display the  
shape of the prior, the likelihood function, and the posterior probability
density function (pdf),
and compare it to the truth.  We see that the mean and mode of the posterior
---both are often taken as single point estimates---when compared to the prior
have moved in the direction of the truth; and the variance of the posterior
(signifying the uncertainty about the true value) is less than that of the
prior.  Another frequent single point estimate is the maximum of the likelihood
function---the well known maximum likelihood estimate $L_{\max}(\kappa)$---which
may be seen as also being close to the truth.  One may also observe that
these (and other) point estimates give only an incomplete picture,
as they will not contain information on the residual uncertainty.  
The Bayesian identification procedure on the other hand yields a complete
probability distribution which informs also about the residual uncertainty.

\begin{figure}[h!]
\centering
\includegraphics*[width=0.68\textwidth]{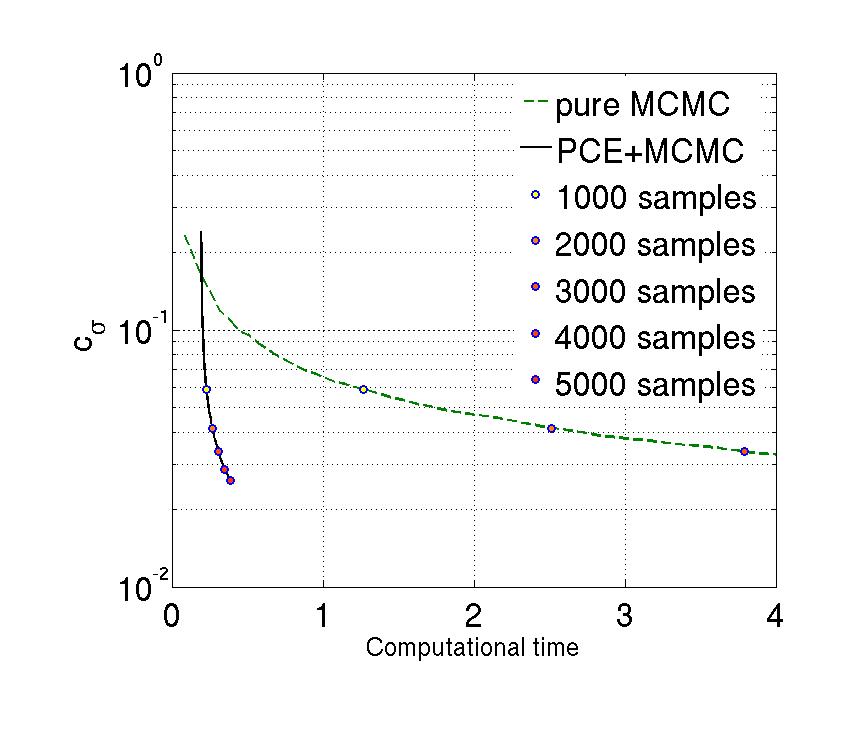}
\caption{Comparison of the computational times of pure MCMC and PCE based MCMC methods versus the estimate $c_\sigma$.}
\label{F:evolstd}
\end{figure}
The MCMC method performs very well as shown, it is computationally 
very expensive.  It requires the 
calculation of the model response for a large number of samples, each time
solving a FE-system with a different material parameter. 
In order to improve the performance \cite{Marzouk2007,Kucherova10}
we compute a polynomial chaos approximation of the model response by a 
stochastic Galerkin procedure.  Then the 
integration of the posterior distribution may be realised by
sampling the PCE as alluded to at the end of \refS{forward}.
This significantly accelerates the update procedure as shown in 
\refig{evolstd}, where the estimate $c_\sigma:=\sigma (L_{\max}(\kappa))
/\EXP{L_{\max}(\kappa)} )$ is plotted versus computational time
for both the pure and the PCE based MCMC approach.  
Clearly the PCE approximation for the same accuracy of estimate runs 
much quicker than the pure MCMC procedure.

\begin{figure}[h!]
\centering
\includegraphics*[width=0.68\textwidth]{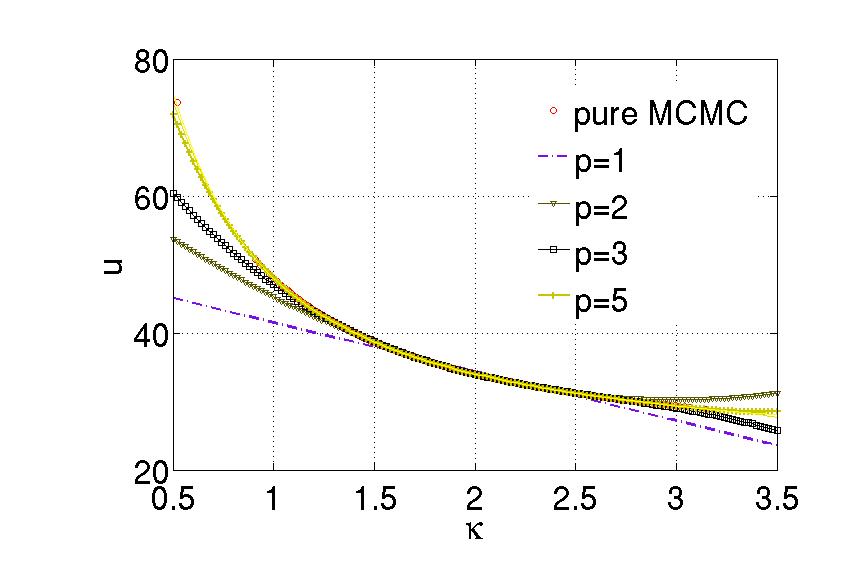}
\caption{Polynomial chaos approximation of the model response in one node of computational domain}
\label{F:pcapprox_1}
\end{figure}
The dependence of the PCE on the polynomial order is investigated in
\refig{pcapprox_1} where the accuracy of the temperature approximation in 
one node of the domain is plotted for different polynomial orders. 
Comparing the approximated solution 
with the reference one obtained by MCMC with $100000$ samples
one concludes that only the fifth order approximation can be safely used. 
In this case both responses match. 
Otherwise, for lower values of 
the polynomial order the mean value is estimated 
correctly (the crossing point) but not the higher order moments 
(the lines do not match).
\begin{figure}[h!]
\centering
\includegraphics*[width=0.68\textwidth]{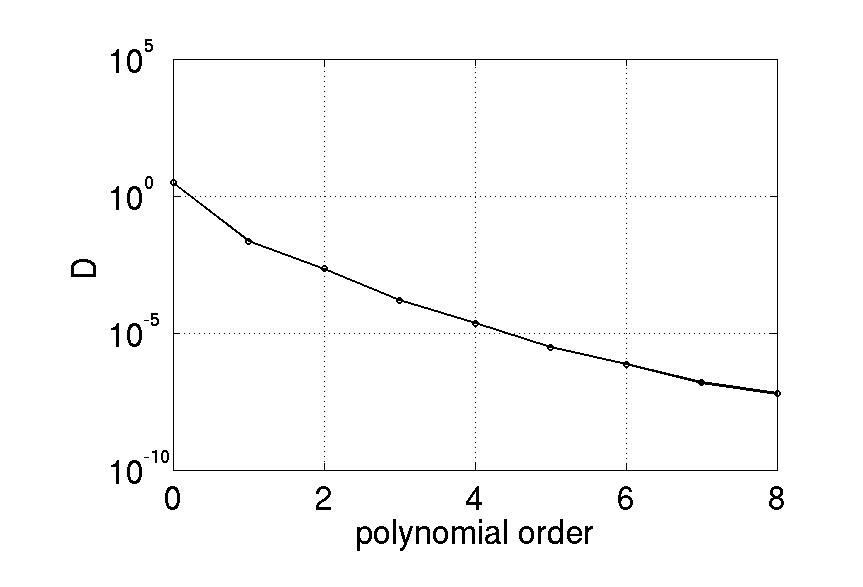}
\caption{Kullback-Leibler distance $D(\tilde{\pi}\|\pi)$ from  the posterior obtained by pure MCMC to the PCE reformulated posterior, versus polynomial order}
\label{F:pcapprox_2}
\end{figure}
We try to quantify this discrepancy in \refig{pcapprox_2} in a 
quantitative assessment of the error in the posterior density of the
PCE-based MCMC by computing the Kullback-Leibler divergence (KLD) of
$\tilde{\pi}_{\kappa}$(posterior obtained by PCE/MCMC) from $\pi_{\kappa}$ (posterior obtained by pure MCMC),
\begin{equation}  \label{eq:KLD}
 D(\tilde{\pi}\|\pi) := \int \tilde{\pi}(\kappa)\log 
   \frac{\tilde{\pi}(\kappa)}{\pi(\kappa)}\,\di\kappa.
\end{equation}
This is shown in \refig{pcapprox_2}, where one may observe a rapid rate
of the convergence of the surrogate posterior $\tilde{\pi}_{\kappa}$ to
the true posterior.

\subsubsection{Conductivity as random field}\label{sec:cgrf}
\begin{figure}[h!]
\centering
\begin{tabular}{cc}
\includegraphics*[width=0.48\textwidth]{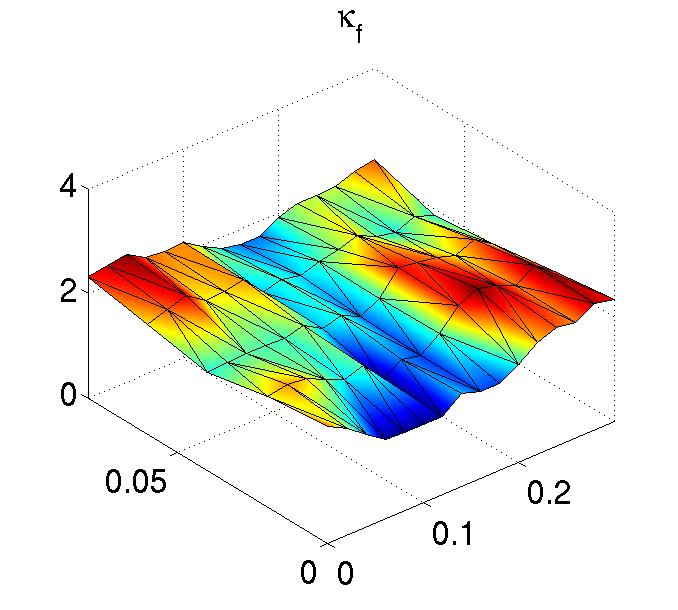}
&
\includegraphics*[width=0.48\textwidth]{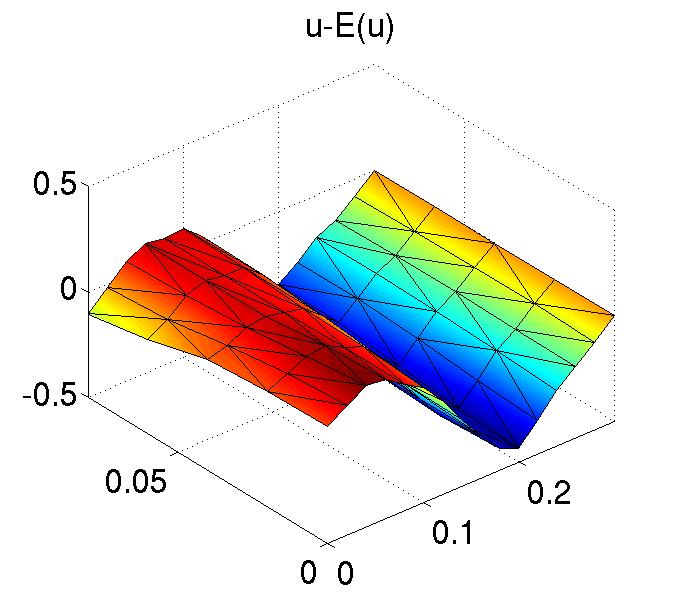}
\\
(a) & (b) 
\end{tabular}
\caption{(a) A realisation of the prior random field, 
  (b) Fluctuations of a corresponding temperature field.}
\label{F:truthfields}
\end{figure}

The next example is a bit more realistic, by not assuming a priori that
the truth is spatially constant.  
We repeat the previous analysis for the prior $\kappa_f$  being a random
field---again for the sake of simplicity---assumed as Gaussian.
The field $\kappa_f=q_f$ is described by a covariance kernel 
$\exp(-r/l_c)$ with $r$ being the distance and 
$l_c$ a correlation length.
Furthermore, the mean and standard deviation are chosen in the same manner 
as before (see Section~\ref{sec:condrv}), while the truth is modelled in 
a more realistic way as one realisation of 
the prior field, as shown in \refig{truthfields}. 
All other quantities such as measurement and measurement error stay the 
same as in the previous example.

\begin{figure}[h!]
\centering
\begin{tabular}{cc}
\includegraphics*[width=0.48\textwidth]{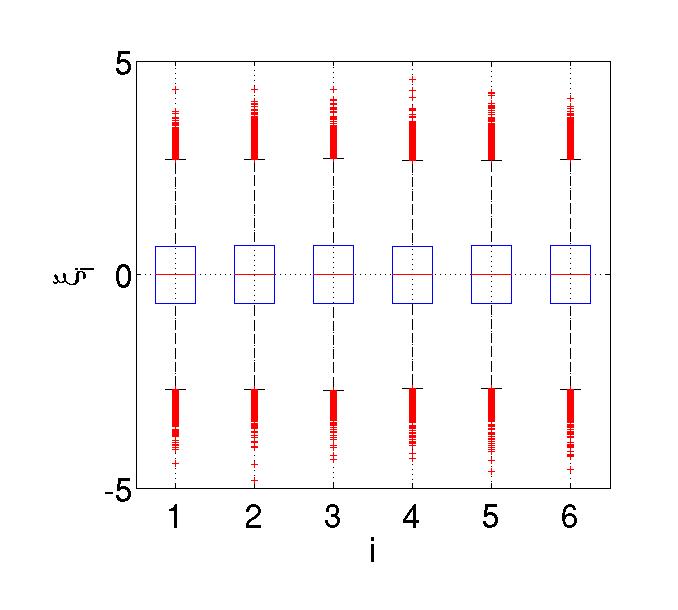} &
\includegraphics*[width=0.48\textwidth]{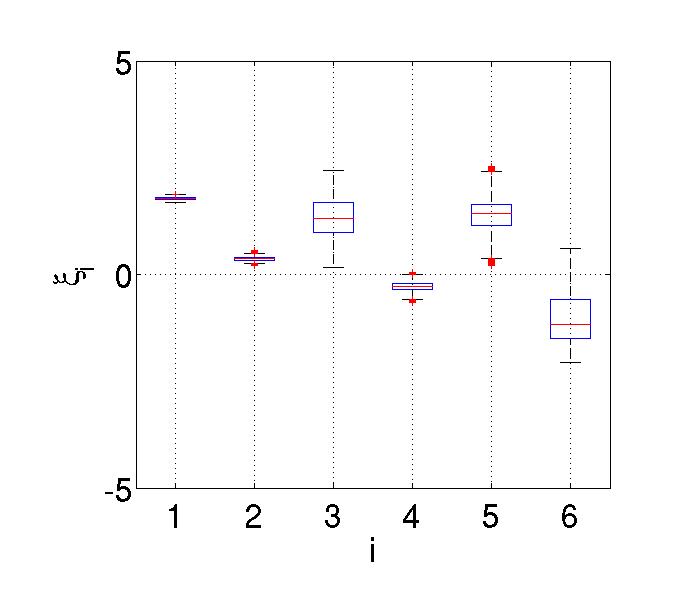}\\
a) & b)\\
\includegraphics*[width=0.48\textwidth]{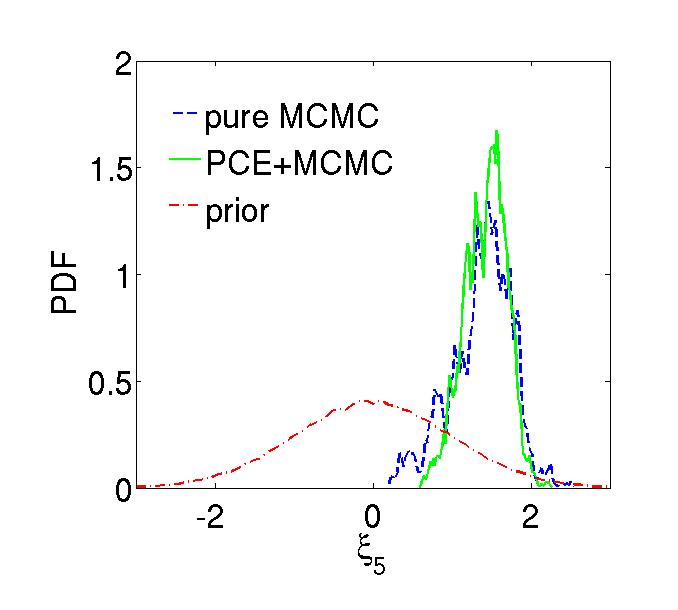} &
\includegraphics*[width=0.48\textwidth]{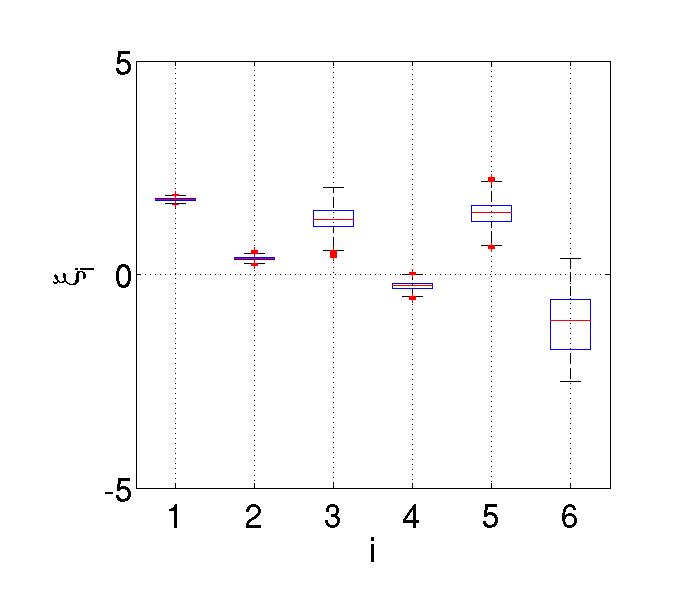}\\
c) &d)\\
\end{tabular}
\caption{Distributions of particular random variables: (a) prior, (b) posterior obtained by pure MCMC, (c) pdfs of $\xi_5$, 
(d) posterior obtained using PCE approximation of model response.  In a), b),
and c), the box graphs show for each RV ($\xi_1,\dots,\xi_6$) on the ordinate
the median, the central 25th percentile as a box, the central 75th
percentile, and outliers.}
\label{F:xi_barpdf}
\end{figure}

The numerical simulation of the conductivity field is performed by truncating 
the prior KLE (see beginning of this \refS{num-xmpls}) 
to six modes, i.e.\  the random variables $\xi_1 \dots \xi_6$ 
in \refeq{eq:exp-KL} have to be updated.  Similarly as before we use the 
pure MCMC update procedure as well as the PCE based MCMC method. 
The pure MCMC Bayesian update uses $100000$ samples in order to
assimilate the posterior variables, while the PCE approximation is taken 
to be of second order. 
In both cases the updated variables obtain the similar non-Gaussian form as 
shown in \refig{xi_barpdf}~b) and d). In addition, by the comparison 
of the bar-graphs of prior and posterior distributions in \refig{xi_barpdf}~b)
and d) one may see the significant reduction of the prior variance. 
This is clearly visible in  \refig{xi_barpdf}~c) where the probability 
density function of $\xi_5$ is plotted. 
Through the identificaton or assimilation process one obtains from a 
very wide prior distribution a much narrower posterior as expected. 
Also, both methods give very similar results leading to the conclusion that
the error introduced by PCE truncation to second order seems to be 
negligible in this
example.  

\subsubsection{Non-linear diffusion} \label{sec:nlindiff}
In the isotropic quasi-static diffusion equation given by \refeq{eq:BAL1}, 
we now assume that the conductivity $\kappa$ depends on the 
temperature field $u(x,\omega)$ through
\begin{equation}
\kappa(u,x,\omega) = \kappa_0(x,\omega) + \kappa_1 u(x,\omega) \, ,
\end{equation}
and the previously linear SPDE becomes nonlinear.
Again, just like the prior assumption of a Gaussian conductivity 
before---chosen purely for reasons of simplicity and useful solely because 
the coefficient of variation was very small---is not consistent 
with the requirement by the second law of thermodynamics
that the conductivity be positive, also this linear dependence of
the conductivity can only be seen as a crude approximation, useful only
for a small temperature range.  Its only purpose is to show how the
method works for a simple nonlinear state equation.

\begin{figure}[htpb]
\centering
\begin{tabular}{cc}
\includegraphics*[keepaspectratio,width=5.5cm]{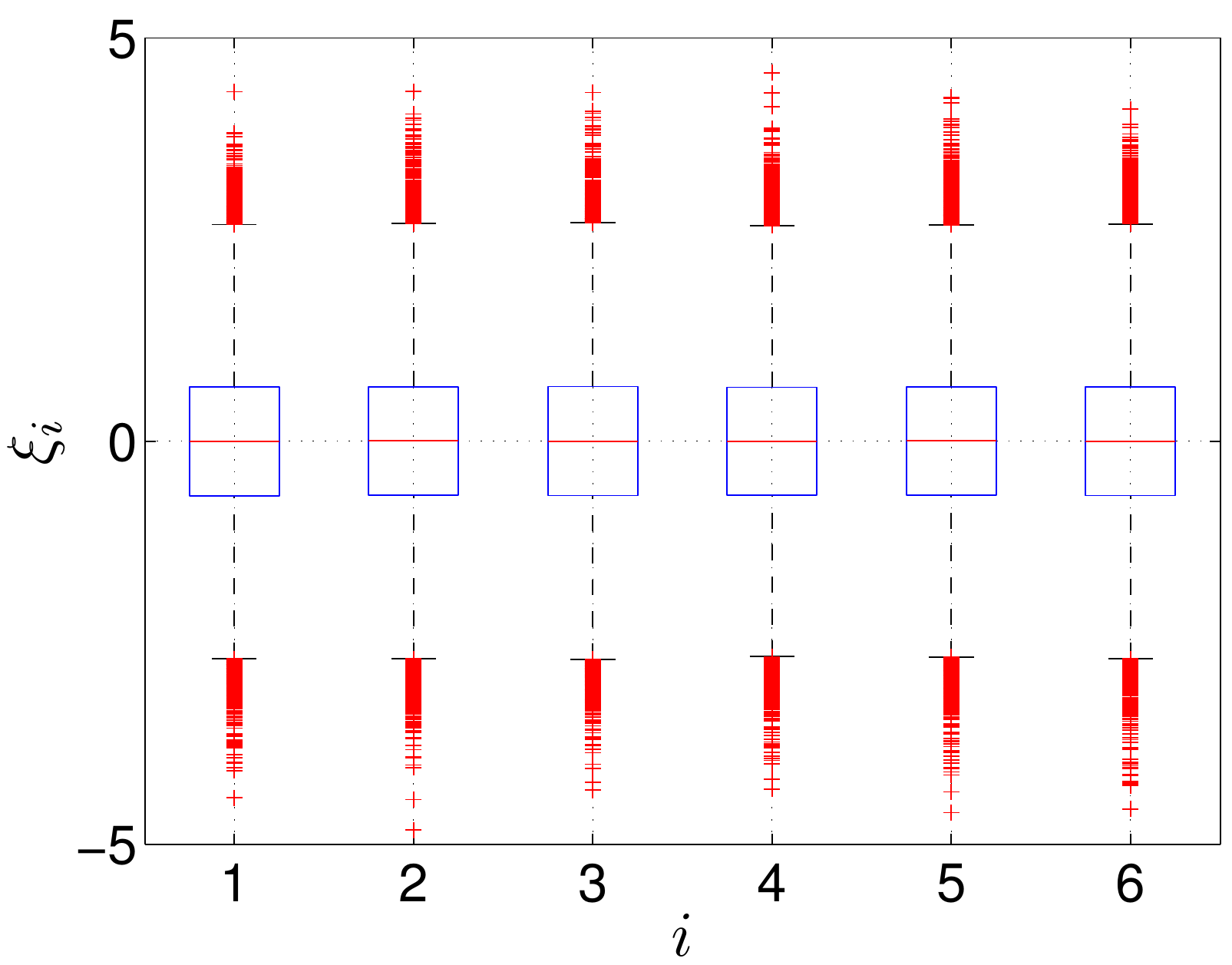}
&
\includegraphics*[keepaspectratio,width=5.5cm]{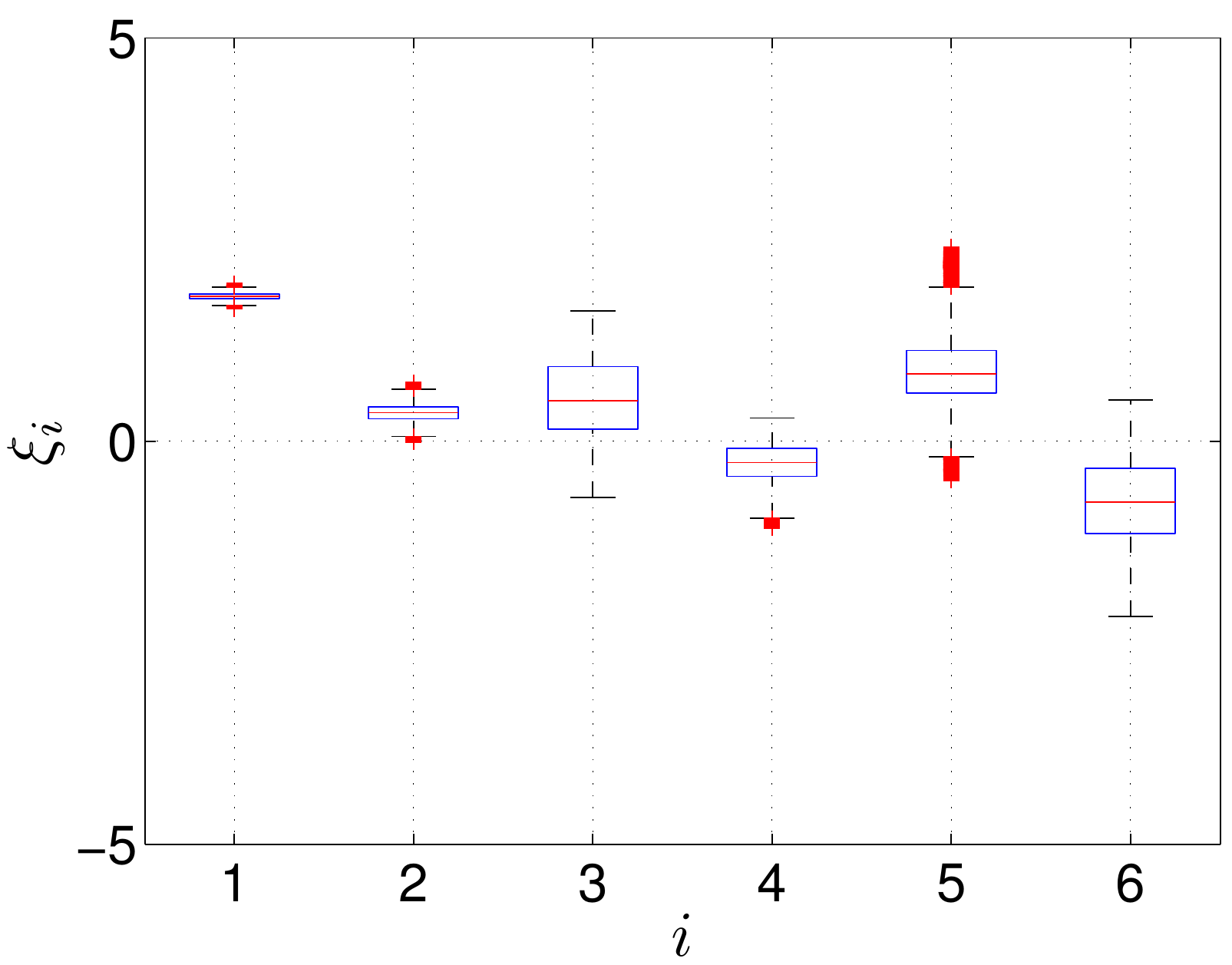}
\\
(a) & (b) \\
\includegraphics*[keepaspectratio,width=5.5cm]{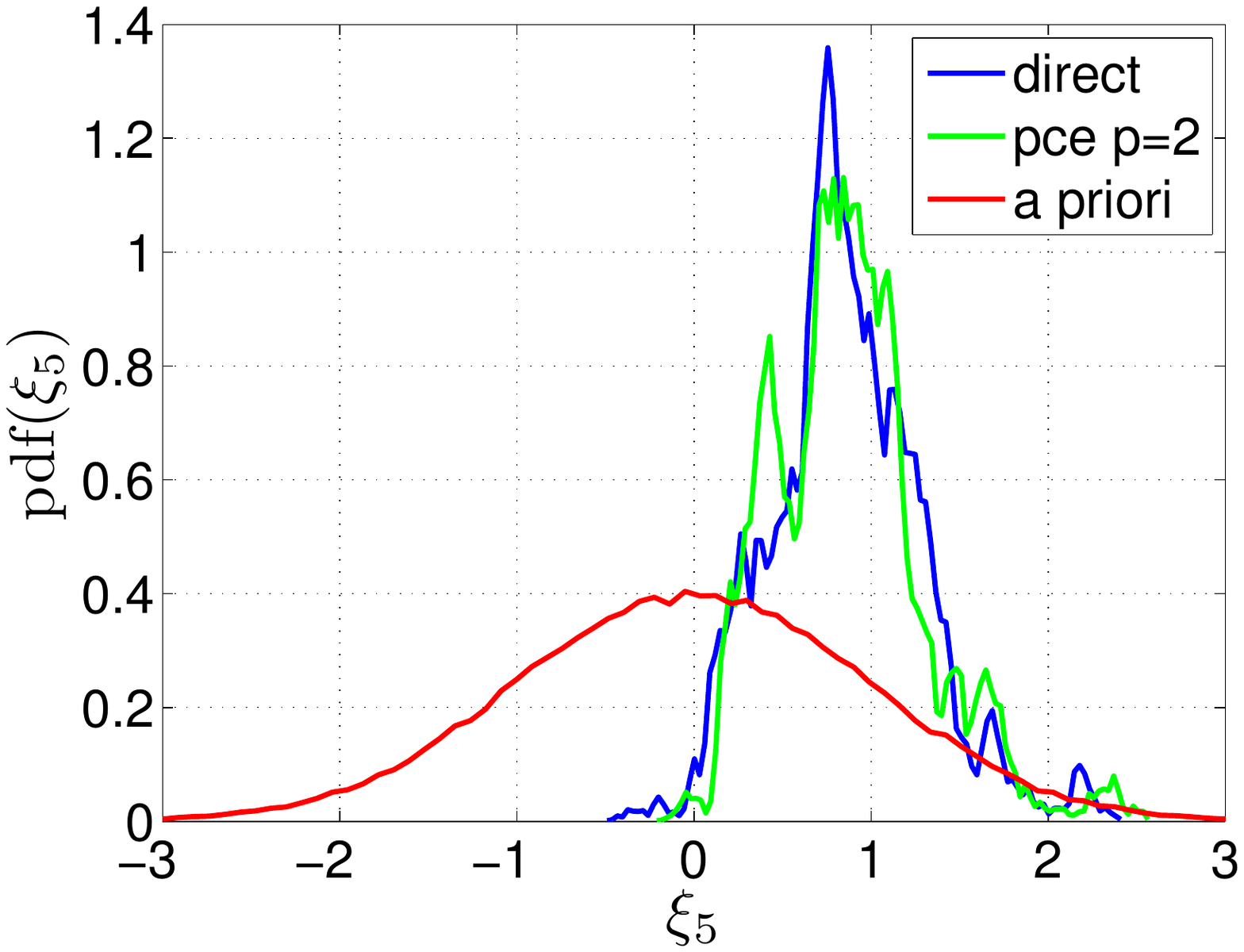}
&
\includegraphics*[keepaspectratio,width=5.5cm]{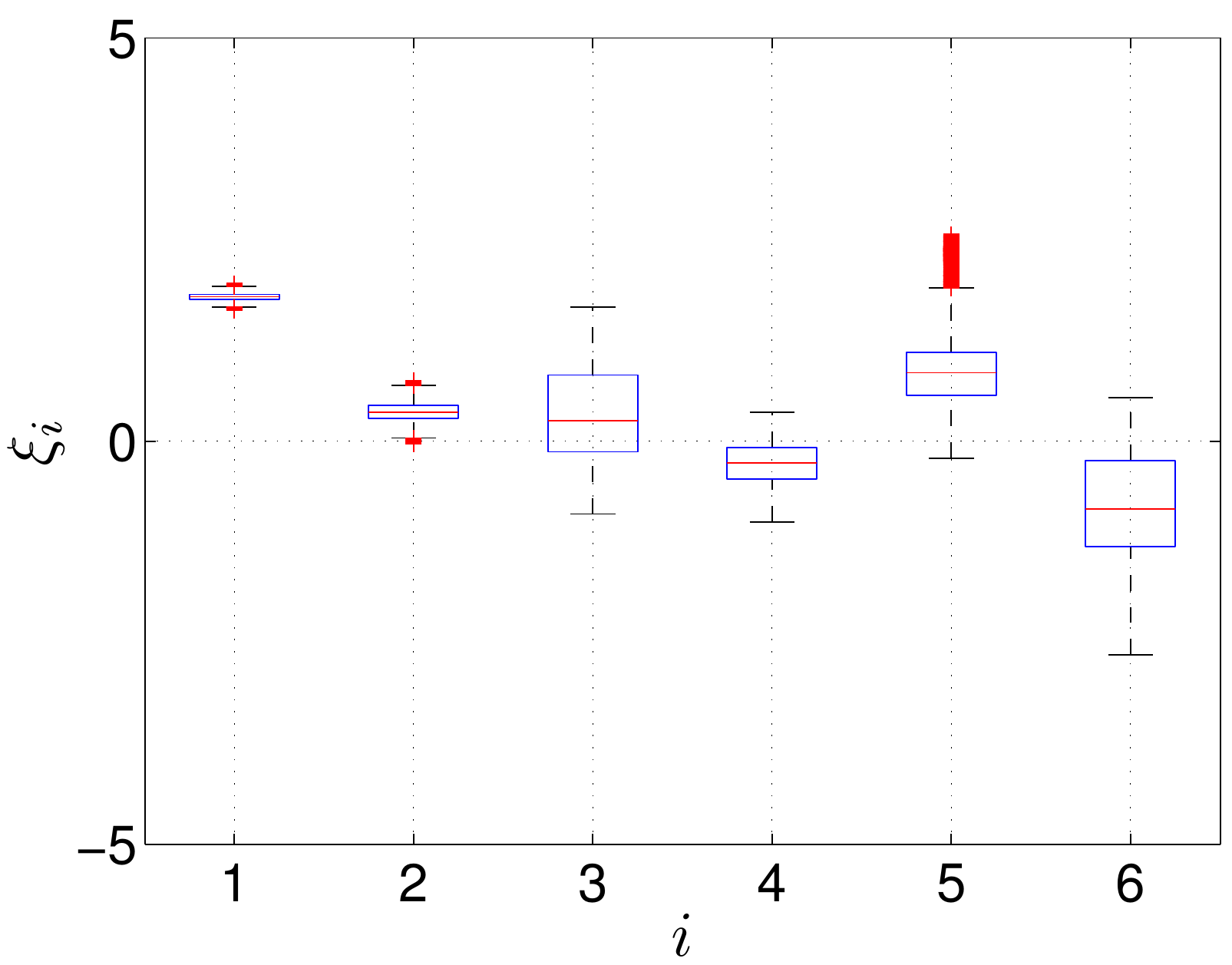}
\\
(c) & (d) 
\end{tabular}
\caption{Distributions of particular random variables: (a) a priori, (b) a posteriori obtained by direct evaluation of the FE model, (c) probability density functions of $\xi_5$, (d) a posteriori obtained using PC approximation of model response.  Box-graphs as in \refig{xi_barpdf}.}
\label{F:nonlpdfs}
\end{figure}

This very simple non-linear constitutive law
is described by two artificial parameters $\kappa_0$ and $\kappa_1$, 
where $\kappa_1 = 0.05\text{ Wm}^{-1}\text{K}^{-2}$ is assumed 
as known and hence is certain,
and $\kappa_0$ is a priori taken as a Gaussian random field
as in the previous Section~\ref{sec:cgrf} with a six term KLE
with mean $ \EXP{\kappa_0}= 2\text{ Wm}^{-1}\text{K}^{-1}$ 
and standard deviation $0.3\text{ Wm}^{-1}\text{K}^{-1}$.  The truth is
taken as one sample of this random field. 
The rest of the experimental set up is taken from the previous tests 
(see the previous Section~\ref{sec:condrv} and \ref{sec:cgrf}).

Similarly as before, the MCMC method in its pure sampling and PCE based form 
is used to calculate the pdfs of the conductivity parameter $\kappa_1$. 
The update procedure consists of assimilation of the six random variables 
$\xi_1,\dots, \xi_6$ of the truncated KLE of the prior field 
$q_f=\kappa_f$  with the help of the temperature measurements in all points 
of the finite element mesh.  Similarly as in the linear case, the box-graphs 
of the prior and posterior random variables in \refig{nonlpdfs}~a) 
and b) show significant reduction of the prior variance and a change 
of the mean, which causes the non-Gaussian form 
of the posterior random variables. 
Furthermore, the comparison of the pure and PCE based MCMC procedure in 
\refig{nonlpdfs}~c) and d) shows that both methods give approximately 
the same results, even though 
the PCE approximation is only of second order.

We see that the MCMC based Bayesian update works also for random fields,
even though the measurement (here the temperature) is not linear in the
quantity to be identified (the conductivity), and also works in case the
state equation is nonlinear.  The inconsistent shape of the prior distribution
(allowing negative conductivities with non-zero probability) is corrected
in all the updates to a form where practically all the probability weight
is in the positive range.  The low order truncated PCE based MCMC performs
almost just as well and is computationally much cheaper---this is 
practically a surrogate model.

\subsection{Examples with linear Bayesian updating} \label{SS:lbudiff}
The examples shown up to now have used the so called `full' Bayesian
update based on \refeq{eq:iII} or \refeq{eq:iIIa}, i.e.\  by updating
the measure.  We now turn our attention to updates which change the
random variable, based on \refeq{eq:proj-t}.  The first example
will also be a diffusion equation, although on a different domain, and
the second example will be an elastic-plastic system, hence non-smooth
and strongly nonlinear.

In previous three examples we have considered the conductivity a priori
as a normally distributed random field, which ignores the positivity requirement
as already noted.  In order to have an a priori model which takes care
of this requirement in a way that our variables are free of 
constraints---the importance of this was already stressed in 
\refSS{cond-exp}---we model the logarithm of the conductivity field; the
a priori field is taken as a lognormal random field, according to maximum 
entropy considerations this is the right choice if one only presumes
that the field has a finite variance.  This makes the identification
problem harder, as the logarithm resp.\ exponential is additionally involved.

As already pointed out in \refSS{cond-exp}, this corresponds with 
\cite{Arsigny_SPDM} considering the positive cone in the space of all
fields as a differentiable Riemann manifold, and more specifically also
as a Lie group.  The logarithm and exponential are then the equally named
maps from Lie theory, carrying the tangent space at the identity into the
corresponding Lie algebra and vice versa.  As the manifold can be equipped
with a Riemann metric which is carried from the Lie algebra via the exponential
map to the tangent space of the identity, and from there via the group
operations to any other place, distances on the manifold can be measured
as path lengths along geodesics.  This turns the manifold into a metric
space and would allow to use the notions of Fr\'echet-type conditional
expectations as alluded to in \refSS{cond-exp}.
But geodesics correspond uniquely to
straight lines in the Lie algebra, hence to their direction vectors without
any constraints \cite{Arsigny_SPDM}.  
It is in this space that we propose to do all the operations.

\subsubsection{Diffusion with linear Bayes updates} \label{SSS:diff-linB}
We consider the same \refeq{eq:BAL1} as before, but now on an $L$-shaped 
domain $\C{G}$ with specified homogenous Dirichlet and Neumann boundary
conditions.  This makes the heat flow a bit more complicated than in the
previous example, where inhomogeneities in the heat flow arise solely due
conductivity variations, whereas here they are also caused through the
geometry.  The external loading is defined as a sinusoidal function
$f=f_0 \sin (\frac{2 \pi}{\lambda} x^T v +\vphi)$,
with $f_0$ being the amplitude, $\lambda$ the wave-length, $\vphi \in [0,2\pi]$
the phase, and $v=(\cos \alpha, \sin \alpha)$  the direction of the
sinusoidal wave specified by an angle $\alpha\in[-\pi/2,\pi/2]$. 
For a more detailed description, please see \cite{bvrAlOpHgm11}.

\begin{figure}[htbp]
\begin{center}
\includegraphics[width=0.68\textwidth]{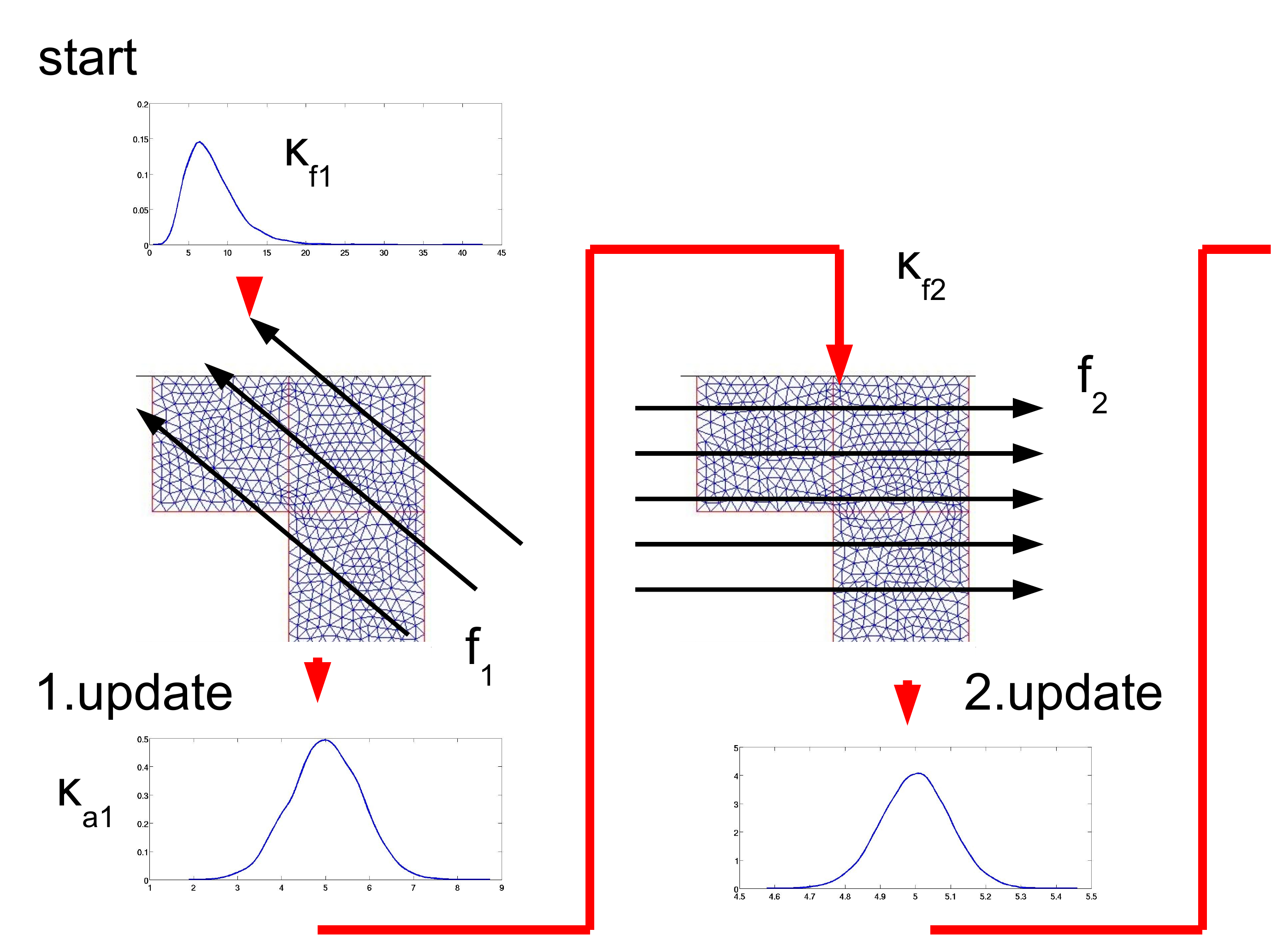}
\caption{Schematic depiction of sequential update procedure}
\label{F:sek_up}
\end{center}
\end{figure}

The identification of the true conductivity field $\kappa_t$ is done with 
the help of the linear Bayesian method in its direct PCE and MC sampling 
(EnKF) form.  Additionally, the update process is realised in a sequential 
way as shown in \refig{sek_up}  by repeating the measurements 
in the same experimental conditions for different values of the right hand side. 
In other words one takes for the prior in the next update the posterior from
the previous update and does the new measurement according to a new value of $f$.
The right hand side values are altered by variation of the appropriate 
parameter values, such as the wave length $\lambda$, the phase $\vphi$, etc. 
In this way different regions of the domain are stimulated, aiding in the
identification process.

\begin{figure}[htbp]
\begin{center}
\begin{tabular}{cc}
\includegraphics[width=0.36\textwidth]{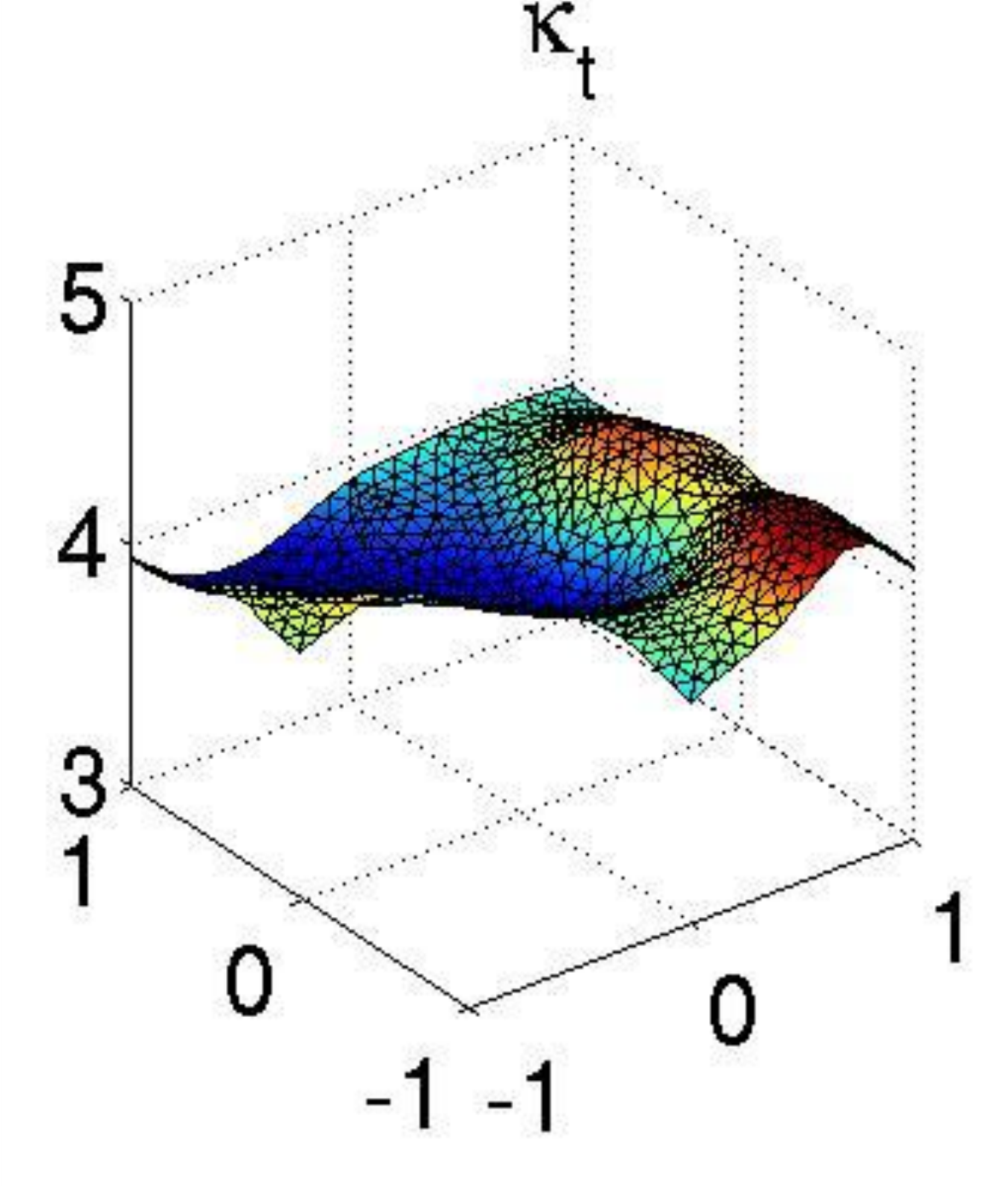} &
\includegraphics[width=0.40\textwidth]{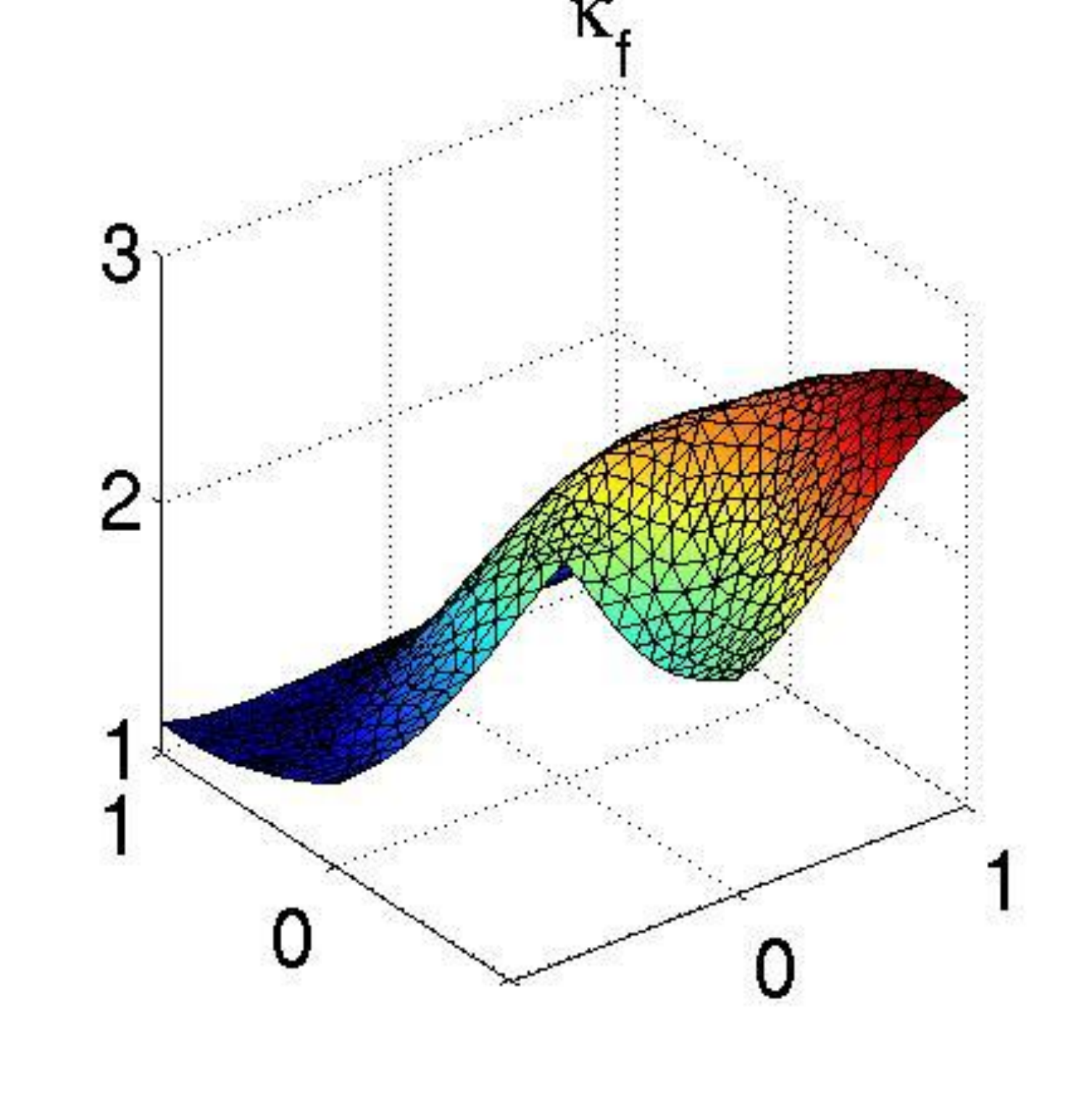}\\
a) & b)\\
\includegraphics[width=0.36\textwidth]{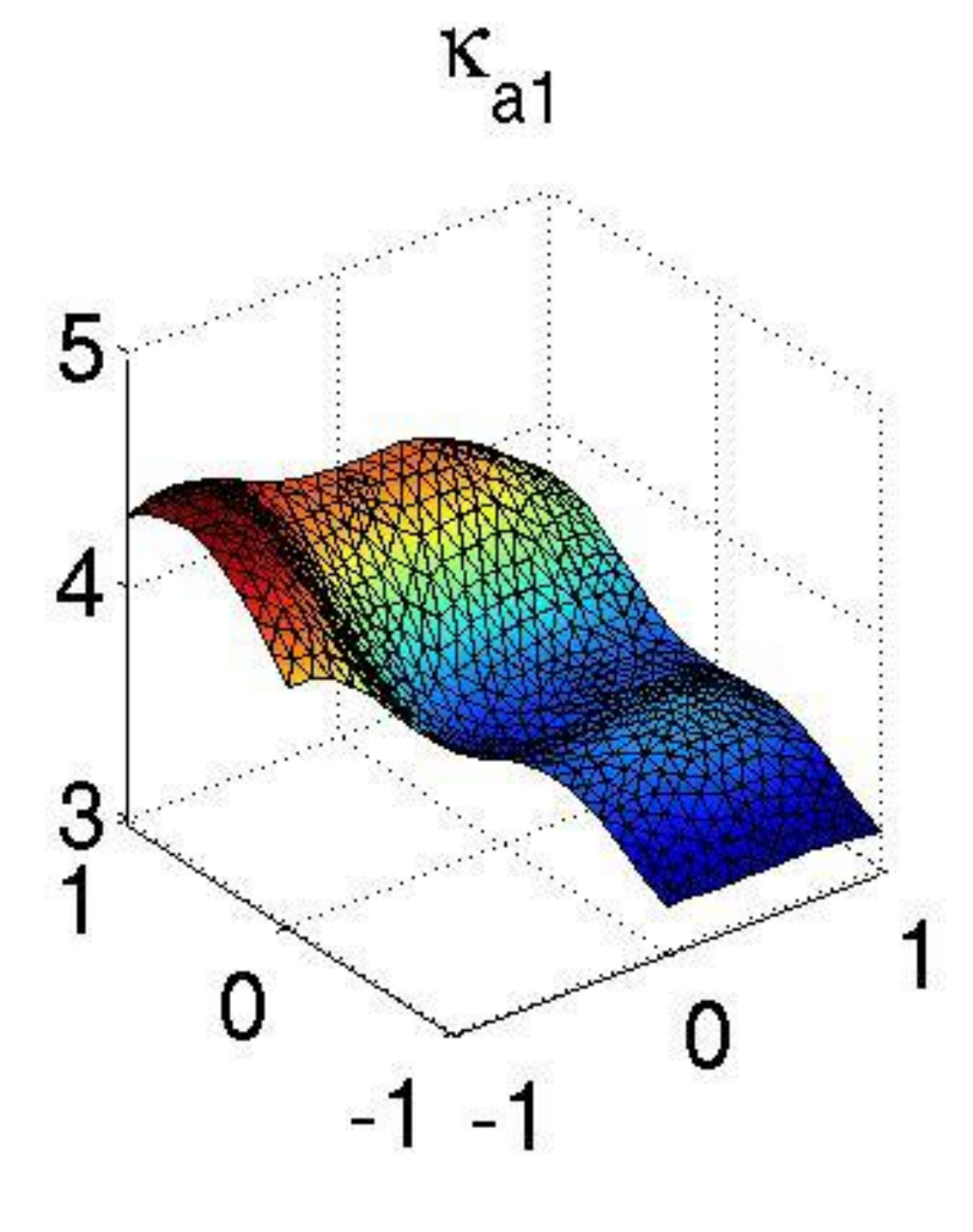}&
\includegraphics[width=0.40\textwidth]{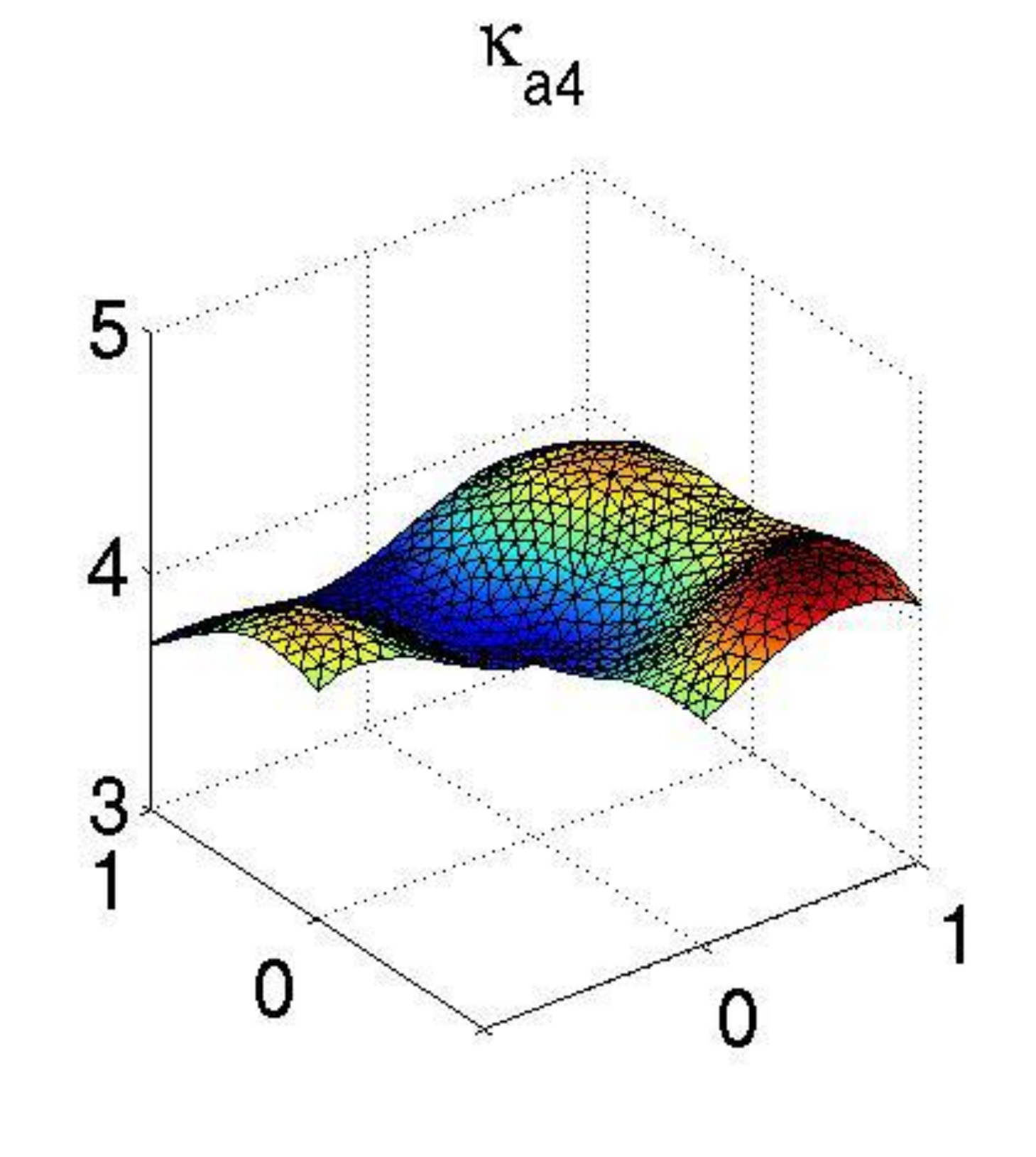}\\
c) & d)\\
\end{tabular}
\caption{Direct PCE update $(p=3,M=10$): a) Truth $\kappa_t$, b) prior, 
c) 1st update $\kappa_{a1}$, d) 4th update $\kappa_{a4}$}
\label{F:truth_with_updates}
\end{center}
\end{figure}

The forward problem is solved within the finite element framework by 
discretising the spatial domain into $1032$ triangular elements. 
The measurements $y(u)$ of the state quantity $u$,
\begin{equation}   \label{eq:bo_measure}
\quad y(u,\omega):=\left[...,y(x_j),...\right]\in \D{R}^{L},\quad 
  y(x_j)=\int_{\C{G}_j} u(x,\omega)\, \di x,
\end{equation}
 are obtained on little patches $\C{G}_j\subset \C{G}$ centred at the 
 finite element node $x_j\in \hat{\C{G}}=\{x_1,...,x_B\}$.
The measurements are performed in only $10\%$ of the total number 
of mesh nodes, equally distributed over the computational area.  The measured
values are disturbed by a centred Gaussian noise with the diagonal covariance 
$\sigma_{\vepsilon}^2\vek{I}$, where $\sigma_\vepsilon$ is approximately equal
to $1\%$ of the measured value.

For the `truth' we take one realisation of a lognormal random field 
sampled from the modified lognormal distribution 
$\kappa_t:=2+\kappa_{b}(x,\omega)$, where $\kappa_b$ represents a 
homogenous lognormal random field with $\EXP{\kappa_b}=1$ and standard
deviation $\sigma_{\kappa_b}=0.2236$ obtained as the exponential transformation 
of a Gaussian random field $q_b$ described by an exponential covariance 
function $\exp(-r/l_c)$, with $r$ being the spatial distance and 
$l_c=0.5$ the correlation length.

Similarly to this, the prior is chosen to be a homogeneous lognormal field  
$\kappa_f=\exp(q_f)$ with $\EXP{\kappa_f}=2.4$ and standard deviation 
$\sigma_{\kappa_f}=0.8944$.  Its  underlying Gaussian random field is 
described by a Gaussian covariance function $\exp(-r^2/l_c^2)$ with a 
correlation length $l_c=2$, which significantly distinguishes 
$\kappa_f$ from $\kappa_t$, see \refig{truth_with_updates}~a) and b).
For the PCE computations we used polynomials with total order up to $p=3$
and $M=10$ KLE modes, while for the EnKF we chose $290$ ensemble members,
as this has approxmately the same workload as the PCE.

\begin{table}[htbp]
        \centering
\begin{tabular}{|c|cc c| c c c|}
\hline
Update & & PCE & &   &EnKF &\\
& $\varepsilon_m$ &  $\bar{\varepsilon}_a$ &$\textrm{var }(\kappa_a)$& $\varepsilon_m$ & $\bar{\varepsilon}_a$ & $\textrm{var }(\kappa_a)$ \\
\hline
A priori & $0.53$ &$0.42$ &$0.8$ &$0.53$ & $0.42$& $0.8$ \\
 $1$st update &$0.18$& $0.14$&$0.13$&$0.19$ &$0.15$ & $0.08$    \\
$2$nd update &$0.07$ &$0.06$&$0.07$ & $0.02$&$0.02$& $0.001$   \\
$3$rd update &$0.06$ &$0.04$&$0.07$ &$0.02$ &$0.02$&  $5.38$e-$04$  \\
$4$th update &$0.03$&$0.03$ &$0.07$ &$0.02$ &$0.02$ &$3.58$e-$04$ \\
\hline
\end{tabular}
\caption{Comparison of relative error of the mode ($\vepsilon_m$), relative
error of the mean ($\bar{\vepsilon}_a$), and the variance 
($\text{var}(\kappa_a)$) of the pdf in one point of the domain 
obtained by the direct (PCE) and the ensemble Kalman filter (EnKF) method.}
\label{reduction}
\end{table}

In  \refig{truth_with_updates}~c) and d) the true field is compared with 
the first and the fourth sequential posterior update obtained by the 
PCE based method.
Since the prior is very different from the truth by all its properties
the `shapes' of the first update and the true field are still relatively 
distinct.  However, the mean and the variance of the prior 
are clearly moved in the direction of the truth.  Doing further updates
with different loadings, the posterior field approaches the truth although
some residual uncertainty remains.  This behaviour is summarised in
\reftab{reduction}, where the PCE and EnKF results are shown for one point 
in the domain. Here one may notice that the EnKF method results in an
unrealistically small variance (smaller than $10^{-2}$  already in the second
update), such that the truth is in a very low probability region.
This undesirable behaviour, which gives such overconfident estimates of
the residual uncertainty is one of the problems of the EnKF, and one has 
to take care not to be misled. 
\begin{figure}[htbp]
\begin{center}
\includegraphics[width=0.69\textwidth]{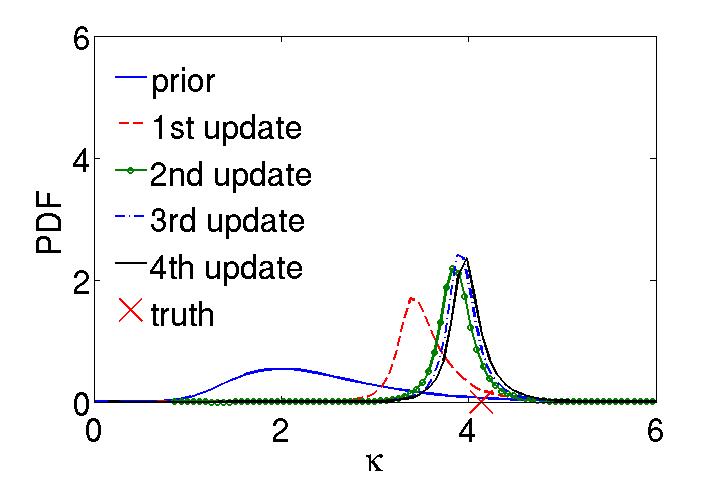}
\caption{Probability density function (pdf) of the posterior for four 
updates by the direct (PCE) method.}
\label{F:pdf_rand_with_updates_1}
\end{center}
\end{figure}
\begin{figure}
\begin{center}
\includegraphics[width=0.71\textwidth]{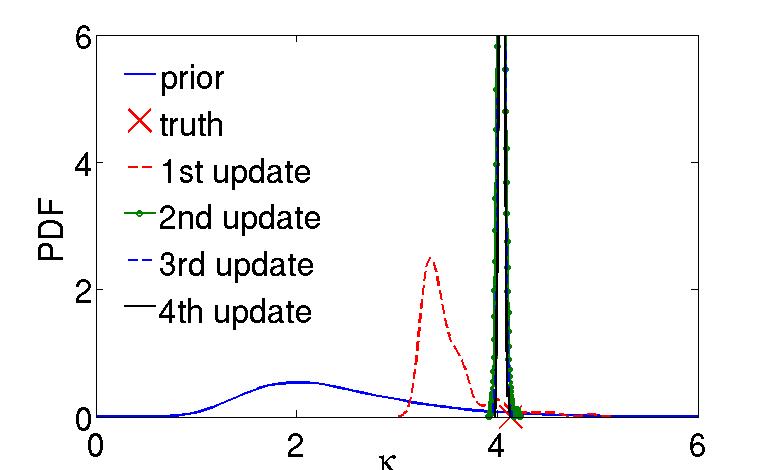}
\caption{Probability density function (pdf) of the posterior for four 
updates by the EnKF method.} 
\label{F:pdf_rand_with_updates_2}
\end{center}
\end{figure}

In a similar vein, comparing the posterior pdfs of the EnKF and the PCE update 
in \refig{pdf_rand_with_updates_1} and \refig{pdf_rand_with_updates_2},
one may observe that after first update the posteriors are fairly similar.
However, in further updates (from second to fourth), the EnKF  appears to 
be very certain about the posterior value, although the truth stays in a very 
low probability region.  In contrast to this, the PCE update is able to keep 
a credible residual variance such
that the truth stays inside a high probability area.

\subsubsection{Elasto-plasticity with linear Bayesian update} \label{SSS:plastic}
To describe quasi-static elasto-plasticity with hardening, we are lead 
a bit beyond the format in \refeq{eq:I} as one has to consider the
non-smooth evolution of the internal parameters 
\cite{HgmRos08a,Rosic_Mat:2008,broHgmMz:eccm10},
which gives a variational inequality as a generalisation of the differential
equation.  
\begin{figure}[htbp]
\begin{center}
\includegraphics[width=0.40\textwidth]{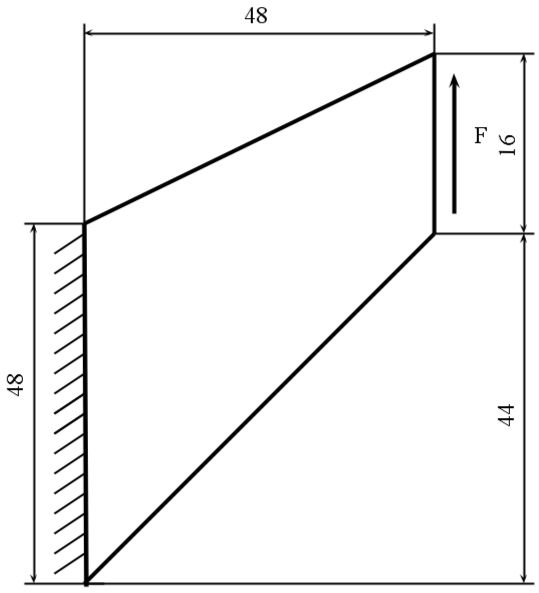}
\caption{Geometry and loading of Cook's membrane test}
\label{F:geometry}
\end{center}
\end{figure}
The state variable is $u = (v, \epsilon_p, \nu)$, where $v$ 
denotes the displacement field, $\epsilon_p$ is the plastic deformation, 
and $\nu$ the appropriate internal hardening variable.  In a mixed formulation
\cite{Rosic_Mat:2008,broHgmMz:eccm10} 
one has to consider at the same time a dual variable
$u^*$ of stress-like quantities.  These quantities have to stay inside
a non-empty closed convex set $\E{K}$, the so-called elastic domain. 
The abstract mixed formulation then is to find functions $u(t)$ and $u^*(t)$
such that $u^*(t) \in \E{K}$ and
\begin{eqnarray} \label{eq:abs_m_1}
 A(u(t))+ u^*(t) = f(t) \nonumber \\
\forall z^*\in \E{K}: \quad \dlangle \dot{u}(t),z^*-u^*(t) \drangle \leq 0.
\end{eqnarray}
In this description, the first equation represents the equilibrium condition,
while the second is the so-called flow rule, stating that the rate of change
$\dot{u}(t) = (\partial/\partial t) u(t)$ is in the normal cone of 
$\E{K}$ at $u^*(t)$; for more details see  
\cite{HgmRos08a,Rosic_Mat:2008,broHgmMz:eccm10}.

\begin{figure}[htbp]
\begin{center}
\begin{tabular}{cc}
\includegraphics[width=0.45\textwidth]{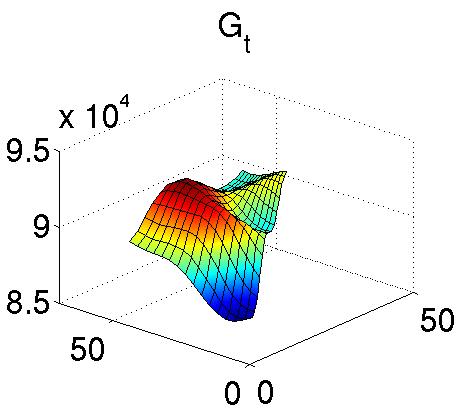} &
\includegraphics[width=0.48\textwidth]{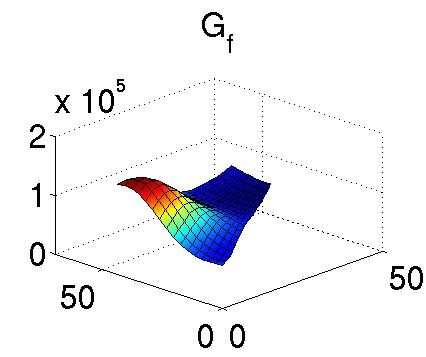} \\
a) & b)\\
\includegraphics[width=0.45\textwidth]{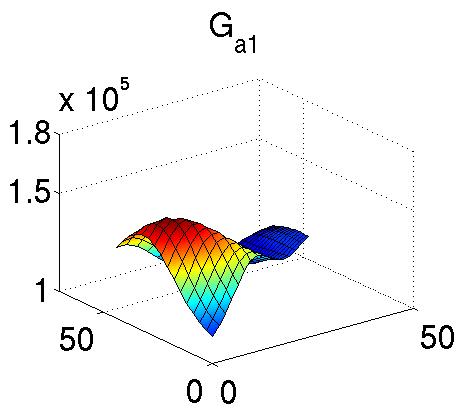} &
\includegraphics[width=0.45\textwidth]{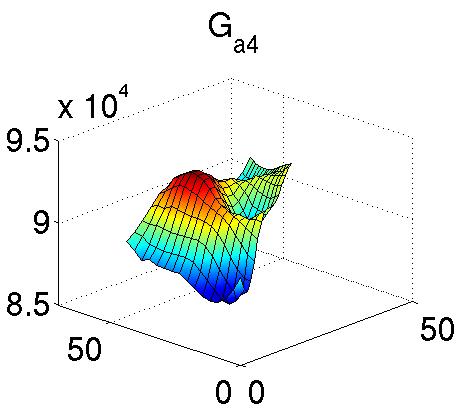}\\
c) & d)\\
\includegraphics[width=0.45\textwidth]{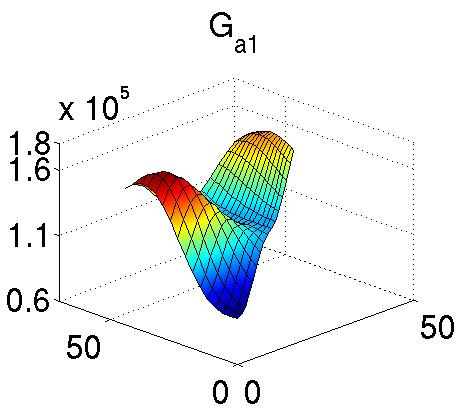} &
\includegraphics[width=0.45\textwidth]{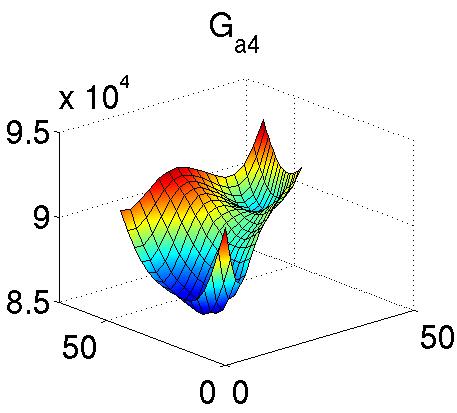}\\
e) & f)\\
\end{tabular}
\caption{Shear modulus $G$ via elastic response: a) truth, b) prior, 
 c) first update, (linear elastic range), d) fourth update (linear elastic range),
 e) first update (nonlinear elasto-plastic range) 
 f) fourth update (nonlinear elasto-plastic range).}
\label{F:elastic_updates}
\end{center}
\end{figure}

Note that while $u^*(t)$ is not on the boundary of $\E{K}$, 
the stress-strain relationship stays linear and the model reduces to 
a linear equation mathematically similar to the diffusion equation,
with $\kappa$ replaced by a constitutive tensor described by the bulk $K$ 
and shear modulus $G$.  Here we investigate the identification of the shear
modulus $G$ via the linear sequential Bayesian procedure 
and shear stress measurements.  The numerical tests are performed on
a well-known and often used example.

Cook's membrane 
is clamped on one end and loaded by a shear force in 
the vertical direction at the other end as shown in \refig{geometry}.
The plate is discretised into $225$ quadrilateral quadratic
eight-noded serendipity finite elements (see \cite{zienkiewiczTaylor00}),
of which $30\%$ are chosen to place the measurements. 
For the sake of simplicity we do not investigate which nodes
are the most appropriate.  Instead, the measurement points are equally 
distributed over the domain, and the measurements are polluted by a
measurement error modelled by independent Gaussian noise with zero mean 
with a diagonal covariance $\sigma_{\varepsilon}^2\vek{I}$,
and $\sigma_{\varepsilon}$ approximately equal to $1\%$ of the measured value.
\begin{table}
        \centering
\begin{tabular}{|c|c c c c c c c|}
\hline
Model &  & & & Update & & &\\ \hline
 & A priori &  $1$st  & $2$nd   & $3$rd  & $4$th & $5$th & $6$th \\ \hline
Elastic & $0.45$  & $0.15$ & $0.04$ & $0.02$ & $0.01$ & $0.01$ & $0.01$\\ \hline
Elasto-Plastic & $0.45$ & $0.30$ & $0.21$ & $0.15$ & $0.08$ & $0.03$ & $0.02$ \\
 \hline
\end{tabular}
\caption{The comparison of the relative root mean square error 
$\vepsilon_a$ in each update for purely elastic and elasto-plastic response.}
\label{elasto_plas_tab}
\end{table}

\begin{figure}[htbp]
\begin{center}
\begin{tabular}{cc}
\includegraphics[width=0.48\textwidth]{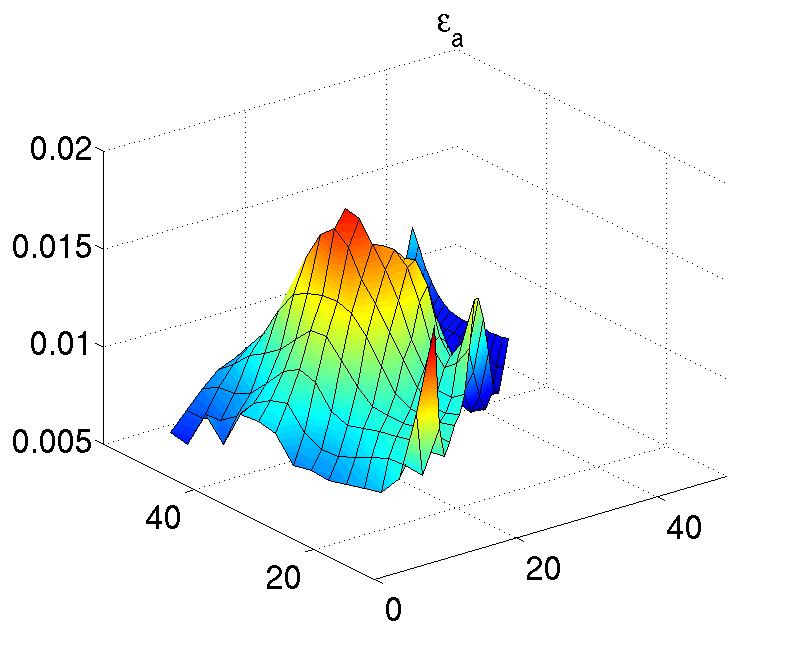}&
\includegraphics[width=0.48\textwidth]{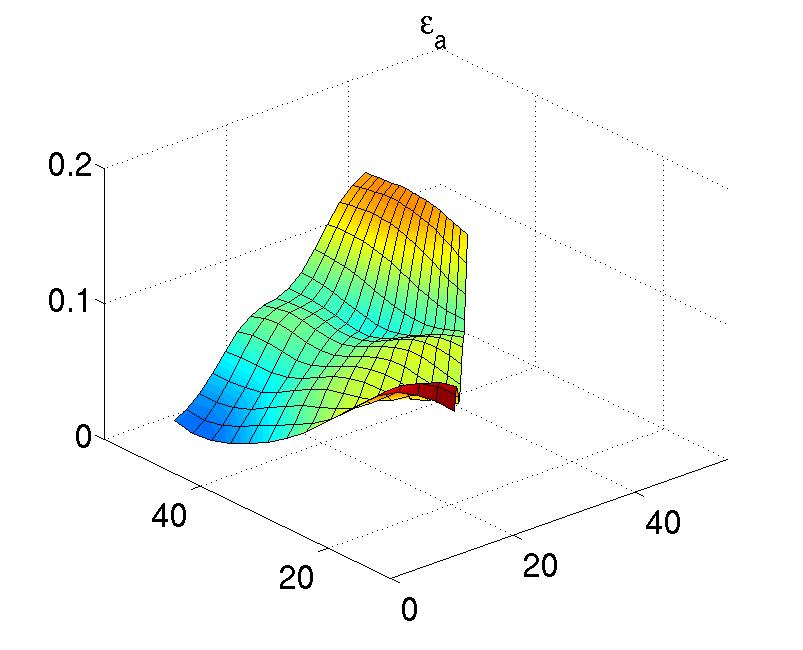}\\
a) &b)\\
\end{tabular}
\caption{Comparison of relative root mean square error $\vepsilon_a$ 
in the fourth update for a) purely elastic, b) elasto-plastic response.}
\label{F:elastic_updates_12}
\end{center}
\end{figure} 

\begin{figure}[htbp]
\begin{center}
\begin{tabular}{c}
\includegraphics[width=0.68\textwidth]{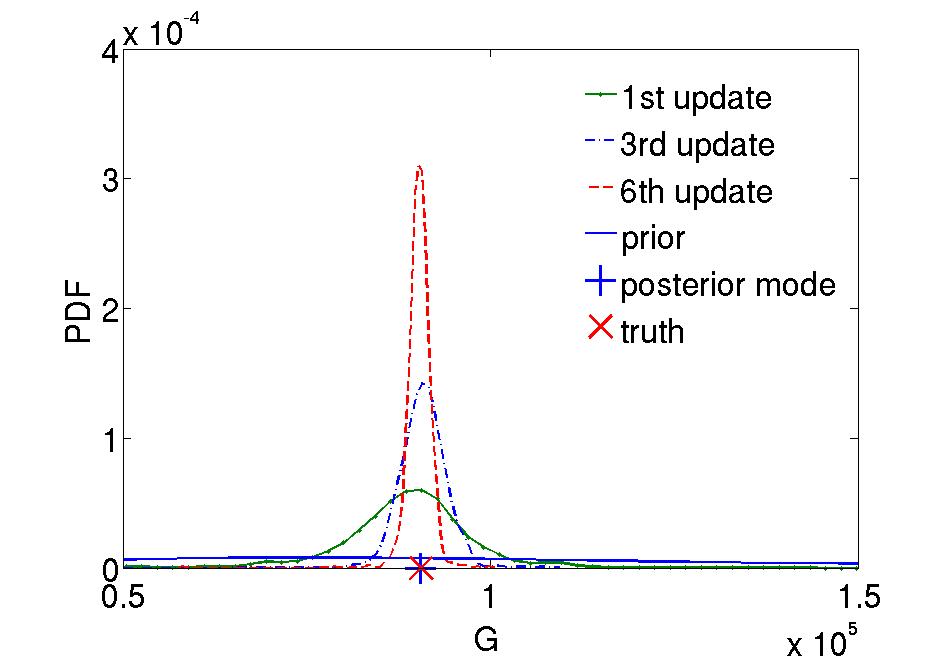}\\
a)\\
\includegraphics[width=0.68\textwidth]{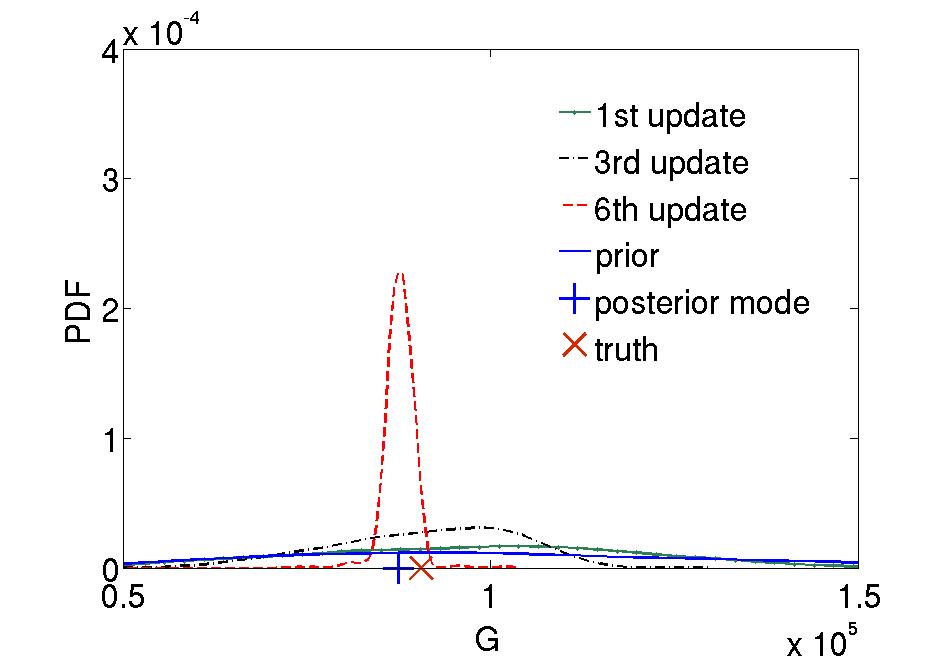}\\
b)\\
\end{tabular}
\caption{Comparison of the posterior and prior pdfs over the updates: 
a) elastic (linear) response, b) elasto-plastic (nonlinear) response.}
\label{F:elastic_updates_12a}
\end{center}
\end{figure}
In this particular example the measurement represents the shear stress $\sigma_{xy}$ numerically computed in a virtual experiment 
with the adopted `true' value of $G_t$.  The truth $G_t$ is taken to 
be one realisation of the modified lognormal field
$G_t=G_0+G_1\kappa_G$, where $\kappa_G:=\exp(q_t)$, $G_0=82760 
\text{ MPa}$, $G_1=0.1 \text{ MPa}$, $\EXP{\kappa_G}=G_0$ 
and standard deviation $\sigma{\kappa_G}=0.1G_0$. 
The underlying Gaussian random field $q_t$ is
described by a Gaussian covariance function $\exp(-r^2/l_c^2)$, 
with $r$ being the spatial distance and $l_c=30$ the correlation length. 

Similarly as for the previous diffusion example in Section~\ref{SSS:diff-linB},
the shear modulus has to be positive and the prior distribution is hence
assumed to be lognormal, though with different characteristics than for the 
truth.  Namely, the prior is assumed to be the lognormal random field 
$G_f=\exp(q_f)$ with $\EXP{G_f}=107590\text{ MPa}$, a standard deviation 
of $0.4 G_0$, and an underlying covariance function of $q_f$
given by $\exp(-r/l_c)$ with $l_c=20$. 
In this way the realisation of the prior is very different from
the realisation of the true field as shown in \refig{elastic_updates}. 

We compute the direct linear Bayesian update via polynomial chaos expansion 
of order $p=3$ and $M=10$ random variables.  The updated field, i.e.\ the
posterior, is then adopted for the new prior in the next sequential update, see
\refig{sek_up}.  Thereby one introduces new information into update process 
by changing  the value of the right hand side. 
With each update the force decreases or increases in every second update 
for $20\%$ starting from $1 \text{ MN/m}$.

The sequential update is done in six stages using both purely elastic and 
elasto-plastic response.  As expected, the nonlinearity is detrimental
to the updating and the posterior obtained via elastic response is better 
than the one from an elasto-plastic model 
as can be seen in \reftab{elasto_plas_tab}, where the relative root mean 
square error reduces faster with each update for purely elastic response.
This is confirmed from the plots of the relative error over the domain.
In \refig{elastic_updates_12} one may see that the elasto-plastic response 
has $10$ times bigger error than the elastic response. In addition,
the error is similar in almost all parts of domain besides the region 
where the plastification starts to occur. 

Finally in \refig{elastic_updates_12a} we compare the pdfs for different
updates for both the elastic/linear and elasto-plastic (nonlinear) response.
The overall picture is similar to what was observed before, namely that
the severe nonlinearity in the elasto-plastic model makes the identification
process much harder.

\section{Conclusion} \label{S:concl}
The problem of identifying parameters or quantities in a computational
model by comparison with either real world measurements or other
computational models (e.g.\ more refined models) is an old and frequent
one.  Practically all approaches start from the idea that the choice of
parameters should be such as to minimise a certain error functional.
Classical methods to do this lead to regularisation.  Another point of
view---taken here---is to embed the unknown quantity in a probability
distribution, where the spread of the probability distribution should
reflect the uncertainty about that quantity.

We have, starting from the elementary textbook formulation of Bayes's
theorem, shown how it connects to conditional expectation.  Conditional
expectation has the minimisation of variance at its background.  This
in turn gives rise to a computational characterisation of updates.
By approximating the space of all measurable functions through the subspace
of linear functions, we disregard a certain amount of information, but
we are rewarded with a simple computation.  Our result contains the well-known
Kalman filter as a special case.

We then discuss how these theoretical constructs can be implemented in
a real computation.  The first group is based on the classical formula
for measures and probability densities.  We sketch the Markov chain Monte
Carlo (MCMC) methods, which may be used in this case.  We also discuss
how the Monte Carlo sampling may be accelerated through the use of
functional approximations for the stochastic part.  We here use standard
Hermite polynomials in Gaussian RVs, i.e.\  Wiener's polynomial chaos.

The second group is based on the conditional expectation idea, leaving the
underlying measure unchanged and updating the relevant RV.  We show how
this idea may be implemented in a sampling Monte Carlo manner---i.e.\ the
ensemble Kalman filter (EnKF)---or alternatively in the PCE setting.  We
then demonstrate the workings of all these methods on some simple examples
of slowly increasing difficulty.  We find that MCMC needs a huge number of
samples, EnKF needs considerably less.  On the other hand the variance
estimate of EnKF seems to be often erring on the optimistic side, thus
severely underestimating the residual uncertainty.  The PCE based methods on
the other hand do not seem to suffer from this illness and give much more
reliable results at a comparable workload.  Our final identification problem
is a very difficult one, an elasto-plastic system, or mathematically speaking
a variational inequality.  The non-smoothness inherent in such problems
slows also the PCE based methods down, but they still succeed.  In our
view these PCE based methods show great promise for these taxing parameter
identification problems.

%
%
%
%




\bibliography{\thebib/jabbrevlong,\thebib/matthies_BU_paper-1,\thebib/phys_D,\thebib/fa,\thebib/risk}




\end{document}